\documentclass{article}
\usepackage[utf8]{inputenc} 
\usepackage{graphicx} 
\usepackage{booktabs} 
\usepackage{array} 
\usepackage{arydshln}
\usepackage{xcolor}
\usepackage{enumerate}
\usepackage{amsmath}
\usepackage{mathtools}
\usepackage{amssymb}
\usepackage{graphicx}
\usepackage{float}
\usepackage{subfig}
\usepackage{fancyhdr}
\usepackage{dcolumn}
\usepackage{tikz,bm,color}
\usetikzlibrary{shapes,arrows}
\usetikzlibrary{calc}
\usepackage[hidelinks]{hyperref}
\usepackage{tikz}
\usetikzlibrary{arrows}
\usetikzlibrary{matrix}
\usetikzlibrary{decorations.pathreplacing}
\usetikzlibrary{calc}
\makeatletter
\newcommand*{\shifttext}[2]{%
  \settowidth{\@tempdima}{#2}%
  \makebox[\@tempdima]{\hspace*{#1}#2}%
}
\newcommand*\dashline{\rotatebox[origin=c]{90}{$\dabar@\dabar@\dabar@$}}
\makeatother
\newcommand{\KW}{\noindent{\bf Keywords: }}
\newcommand{\AMS}{\medskip\noindent{\bf Mathematics Subject Classification 2010: }}

\newcommand{\uparrowbarred}{\mathbin{\rotatebox[origin=c]{90}{\shifttext{8pt}{$\mapstochar$}\shifttext{2pt}{$\rightarrow$}}}}
\newcommand{\downarrowbarred}{\mathbin{\rotatebox[origin=c]{-90}{\shifttext{2pt}{$\mapstochar$}\shifttext{-2pt}{$\rightarrow$}}}}

\makeatletter
\newlength\xvec@height%
\newlength\xvec@depth%
\newlength\xvec@width%
\newcommand{\xvec}[2][]{%
  \ifmmode%
    \settoheight{\xvec@height}{$#2$}%
    \settodepth{\xvec@depth}{$#2$}%
    \settowidth{\xvec@width}{$#2$}%
  \else%
    \settoheight{\xvec@height}{#2}%
    \settodepth{\xvec@depth}{#2}%
    \settowidth{\xvec@width}{#2}%
  \fi%
  \def\xvec@arg{#1}%
  \def\xvec@dd{:}%
  \def\xvec@d{.}%
  \raisebox{.2ex}{\raisebox{\xvec@height}{\rlap{%
    \kern.05em
    \begin{tikzpicture}[scale=1]
    \pgfsetroundcap
    \draw (.05em,0)--(\xvec@width-.05em,0);
    \draw (\xvec@width/2,-0.15em)--(\xvec@width/2, 0.15em);
    \ifx\xvec@arg\xvec@d%
      \fill(\xvec@width*.45,.5ex) circle (.5pt);%
    \else\ifx\xvec@arg\xvec@dd%
      \fill(\xvec@width*.30,.5ex) circle (.5pt);%
      \fill(\xvec@width*.65,.5ex) circle (.5pt);%
    \fi\fi%
    \end{tikzpicture}%
  }}}%
  #2%
}
\newcommand{\xvecsub}[2][]{%
  \ifmmode%
    \settoheight{\xvec@height}{$#2$}%
    \settodepth{\xvec@depth}{$#2$}%
    \settowidth{\xvec@width}{$#2$}%
  \else%
    \settoheight{\xvec@height}{#2}%
    \settodepth{\xvec@depth}{#2}%
    \settowidth{\xvec@width}{#2}%
  \fi%
  \def\xvec@arg{#1}%
  \def\xvec@dd{:}%
  \def\xvec@d{.}%
  \raisebox{-0.3ex}{\raisebox{\xvec@height}{\rlap{%
    \kern.05em
    \begin{tikzpicture}[scale=1]
    \pgfsetroundcap
    \draw (.05em,0)--(\xvec@width-.05em,0);
    \draw (\xvec@width/2,-0.1em)--(\xvec@width/2, 0.1em);
    \ifx\xvec@arg\xvec@d%
      \fill(\xvec@width*.45,.5ex) circle (.5pt);%
    \else\ifx\xvec@arg\xvec@dd%
      \fill(\xvec@width*.30,.5ex) circle (.5pt);%
      \fill(\xvec@width*.65,.5ex) circle (.5pt);%
    \fi\fi%
    \end{tikzpicture}%
  }}}%
  #2%
}
\makeatother
 
 \DeclareRobustCommand{\naturalto}{%
  \mathrel{
  \ooalign{$\longrightarrow$\cr%
  \kern1.45ex\raise.275ex\hbox{\scalebox{1}[0.522]{$\mid$}}\cr}
  }%
}

 \DeclareRobustCommand{\naturaltolr}{%
  \mathrel{
  \ooalign{$\longleftrightarrow$\cr%
  \kern1.75ex\raise.275ex\hbox{\scalebox{1}[0.522]{$\mid$}}\cr}
  }%
}

\begin{document}
\title{Classical Bivalent Logic as a Particular Case of Canonical Fuzzy Logic} \author{
Osvaldo Skliar\thanks{Universidad Nacional, Costa Rica. E-mail: osvaldoskliar@gmail.com. https://orcid.org/0000-0002-8321-3858.}
\and Sherry Gapper\thanks{Universidad Nacional, Costa Rica. E-mail: sherry.gapper.morrow@una.ac.cr. https://orcid.org/0000-0003-4920-6977.}
\and Ricardo E. Monge\thanks{Universidad CENFOTEC, Costa Rica. E-mail: rmonge@ucenfotec.ac.cr. https://orcid.org/0000-0002-4321-5410.}}


\maketitle
\begin{abstract}
A review is presented of the correspondence existing in both classical bivalent logic (BL) and canonical fuzzy logic (CFL) between each law or tautology in propositional calculus and a law in set theory. The latter law consists of the equality of a) a set whose structure is isomorphic to the law considered in propositional calculus and b) the universal set. In addition to the operations of CFL considered previously by the authors, initial attention is given to the operations with infinite sets by considering two of them: the union of sets and the intersection of sets. Attention is also given to how propositional calculus in BL can be considered a particular case of propositional calculus in CFL, and how the theory of classical sets can be considered a particular case of the theory of fuzzy sets according to CFL.
\end{abstract}

\KW classical bivalent logic, canonical fuzzy logic, theory of classical sets, theory of canonical fuzzy sets

\AMS 03B05, 03B52, 03E20, 03E72, 03E99

\newpage
\section{Introduction}
This article has been written as part of an ongoing research program on the following topics:\\

\noindent{I. the possibility of considering classical bivalent logic (BL) as a particular ``limit'' case of certain variants of diverse non-classical logics; and}

\noindent{II. the relations existing between different calculi, such as propositional calculus and predicate calculus, both in BL and in the variants mentioned.}\\

The objectives of this article are the following:\\

\noindent{1) Briefly review previously presented results on Canonical Fuzzy Logic (CFL) \cite{a1} -- a variant of fuzzy logic;
}

\noindent{2) show how, both in BL and CFL, each law (or tautology) in propositional calculus is isomorphic to a law (or tautology) in the respective set theories -- i.e., theory of classical sets and CFL set theory;}

\noindent{3) approach a topic not covered in \cite{a1}: the operations with infinite sets (that is, sets to which infinite elements belong) in fuzzy set theory according to CFL. The operations of union and intersection of sets will be addressed; and}

\noindent{4) specify a) how propositional calculus in BL can be considered a particular case of propositional calculus in CFL, and b) how the theory of classical sets can be considered a particular case of the theory of fuzzy sets in CFL.}

For a clear understanding of this article, only basic knowledge of propositional calculus and of set theory is required. As an introduction to logic, whose basic calculus is propositional calculus, one may consult, for example, \cite{b1}, \cite{b2}, \cite{b3}, \cite{b4} and \cite{b5}. On set theory, one may consult \cite{c1}, \cite{c2}, \cite{c3} and \cite{c4}; and on fuzzy logic and fuzzy sets, \cite{d1}, \cite{d2}, \cite{d3}, and \cite{d4}.

\section{A Correspondence in BL between Each Operation in Propositional Calculus and a Particular Operation in Set Theory}

In BL every proposition, or statement, can have only one of two possible truth values: It can only be true or false. The negation of a true proposition is a false proposition and the negation of a false proposition is a true proposition.

Propositional variables, variables which can be replaced by propositions, are usually denominated $p$, $q$, $r$, $\ldots$. In this article the letter $p$ will not be used to denominate any propositional variable to prevent confusion in future related articles in which different probabilities will be denominated $p_1$, $p_2$, $p_3$, $\ldots$. Reference will be made to the different propositional variables as $q_1$, $q_2$, $q_3$, $\ldots$. In addition, use will be made of a ``license'', often found in literature on logic, to refer to those propositional variables as if they were propositions. Thus, when stating, ``$q_1$ is true'' or ``$q_2$ is false'', for example, one is implicitly stating, respectively, ``Admit that $q_1$ has been replaced by a true proposition'', and ``Admit that $q_2$ has been replaced by a false proposition''.

Two sets have an important role in set theory: the universal set (or the universe of discourse) $\mathbb U$; and the empty set $\varnothing$.

All the elements of $\mathbb U$ that do not belong to $C$, characterized within the frame of that $\mathbb U$, belong to the complement of $C$, which is the set $\xvec{C}$. Each element of that $\mathbb U$ can belong -- or not belong -- to each set considered within a frame of a given $\mathbb U$. Those sets are denominated $C_1$, $C_2$, $C_3$, $\ldots$. If reference is made to only one of those sets, the subscript is unnecessary and it can simply be denominated $C$. If all the elements of $\mathbb U$ belong to a set $C$, that $C$ is equal to $\mathbb U$: $C = \mathbb U$. If no element of $\mathbb U$ belongs to a given set $C$, that $C$ is equal to $\varnothing$: $C = \varnothing$.

If the same elements belong to each of any two sets $C_1$ and $C_2$ whatsoever characterized within the frame of a given $\mathbb U$, then those sets are equal: $C_1 = C_2$.

To negate a proposition $q$, use can be made of the operator of \textit{negation} in propositional calculus, which will be symbolized by a horizontal bar placed over the proposition that it negates. Thus, $\overline{q}$ (or not $q$) is the result of the action of the operator of negation on the operand which is the proposition $q$.

The operator of \textit{complementation} of any set $C$ whatsoever will be represented by the symbol $\xvec{\quad}$ placed over a set whose complement is sought. Therefore, the set which is a complement of $C$ is symbolized as $\xvec{C}$. When the operator of complementation acts on the operand $C$, the result is the complement of $C$, which is the set $\xvec{C}$. All the elements belonging to $\mathbb U$ that do not belong to $C$ belong to the complement of any set $C$, which is the set $\xvec{C}$.

The proposition $\overline{\overline{q}}$ is the negation of the negation, or ``double negation'', of the proposition $q$. The set $\xvec{\xvec{C}}$ is the complement of the complement of $C$.
 
The truth tables for $\overline{q}$ and $\overline{\overline{q}}$ are shown in figure \ref{f1a}, and the membership tables for $\xvec{C}$ and $\xvec{\xvec{C}}$ are shown in figure \ref{f1b}.

 \begin{figure}[H]
\centering
\subfloat[Truth tables for $\overline{q}$ and $\overline{\overline{q}}$]{
\hspace{2in}
\begin{tabular}{c||c||c}
$q$ & $\overline{q}$ & $\overline{\overline{q}}$ \\
\midrule
0 & 1  & 0 \\ 
1 & 0  & 1 \\ 
\end{tabular}
\hspace{2in}
\label{f1a}
}
\\
\subfloat[Membership tables for the sets $\xvec{C}$ and $\xvec{\xvec{C}}$]{
\hspace{1.9in}
\begin{tabular}{c||c||c}
$C$ & $\xvec{C}$ & $\xvec{\xvec{C}}$  \\
\midrule
0 & 1 & 0  \\ 
1 & 0  & 1 \\ 
\end{tabular}
\hspace{2in}
\label{f1b}}%
\caption{a) Truth tables for $\overline{q}$ and $\overline{\overline{q}}$, and b) Membership tables for $\xvec{C}$ and $\xvec{\xvec{C}}$}
\label{f1}
\end{figure}

The presence of an element 0 in any column whatsoever of the truth tables shown in figure \ref{f1a} means that the proposition to which that column corresponds is considered false. The presence of a 1 in any column whatsoever of the truth tables shown in figure \ref{f1a} means that the proposition to which that column corresponds is considered true.

The numerical sequence of the first row of the truth tables in figure \ref{f1a} (0, 1, 0) should be interpreted as follows: If $q$ is true, then $\overline{q}$ is false and $\overline{\overline{q}}$ is true.

Note that the column corresponding to $\overline{\overline{q}}$ is equal to the column corresponding to $q$. This means that the negation of the negation of $q$ -- that is, the ``double negation'' of $q$ -- has the same truth value as $q$. (Actually it is identical to $q$.)

The numerical sequence of the first row of the membership tables shown in figure \ref{f1b} (0, 1, 0) should be interpreted as follows: If any element whatsoever belonging to the universal set $\mathbb U$ does not belong to $C$, as specified by the first number in that sequence (0), then that element does belong to the set $\xvec{C}$, as specified by the second number in that sequence (1), and does not belong to the set $\xvec{\xvec{C}}$, as specified by the third number in that sequence (0).

The numerical sequence of the second row of the membership tables shown in figure \ref{f1b} (1, 0, 1) should be interpreted as follows: If any element whatsoever belonging to the universal set $\mathbb U$ considered belongs to the set $C$, as specified by the first number in that sequence (1), then that element does not belong to the set $\xvec{C}$, as specified by the second number of that sequence (0), and does belong to the set $\xvec{\xvec{C}}$, as specified by the third number in that sequence (1).

Note that the column corresponding to $\xvec{\xvec{C}}$ is equal to the column corresponding to $C$. This means that the complement of the complement of $C$ -- that is, the ``double complement'' of $C$ -- is equal to $C$: $C = \xvec{\xvec{C}}$.

From the viewpoint of the presence of zeros and ones, the truth table for $\overline{q}$ is equal to the membership table for $\xvec{C}$, and the truth table for $\overline{\overline{q}}$ is equal to the membership table for $\xvec{\xvec{C}}$. This result makes it possible to establish the following correspondences: 1) between $q$ and $C$, 2) between the operator of negation in the propositional calculus of BL and the operator of complementation in set theory, 3) between $\overline{q}$ and $\xvec{C}$, and 4) between $\overline{\overline{q}}$ and $\xvec{\xvec{C}}$.

The operator of \textit{disjunction} of the two propositions $q_1$ and $q_2$ is symbolized as $\lor$. The proposition $q_1 \lor q_2$ is read as ``$q_1$ or $q_2$''.

The operator of the \textit{union} of two sets $C_1$ and $C_2$ is symbolized as $\cup$. The set $C_1 \cup C_2$ is the set resulting from the union of the sets $C_1$ and $C_2$.

The truth table for the proposition $q_1 \lor q_2$ is shown in figure \ref{f2a}, and the membership table for the set $C_1 \cup C_2$ is shown in figure \ref{f2b}.

\begin{figure}[H]
\centering
\subfloat[Truth table for $q_1 \lor q_2$.]{
\hspace{1.6in}
\begin{tabular}{c|c||c|}
$q_1$ &  $q_2$ & $q_1 \lor q_2$  \\
\midrule
0 & 0 & 0   \\  
0 & 1 & 1   \\  
1 & 0 & 1   \\  
1 & 1 & 1   \\  
\end{tabular}
\hspace{2in}
\label{f2a}
}
\\
\subfloat[Membership table for the set $C_1 \cup C_2$.]{
\hspace{1.5in}
\begin{tabular}{c|c||c|}
$C_1$ &  $C_2$ & $C_1 \cup C_2$  \\
\midrule
0 & 0 & 0   \\  
0 & 1 & 1   \\  
1 & 0 & 1   \\  
1 & 1 & 1   \\  
\end{tabular}
\hspace{2in}
\label{f2b}}%
\caption{a) Truth table for $q_1 \lor q_2$, and b) Membership table for $C_1 \cup C_2$}
\label{f2}
\end{figure}

In the first row of the truth table in figure \ref{f2a} the numerical sequence shown is 0, 0, 0. It should be interpreted as follows: If $q_1$ is false, and $q_2$ also is false, then the disjunction of $q_1$ and $q_2$ -- that is, the proposition $q_1 \lor q_2$ -- is false as well. Note that the case in which both $q_1$ and $q_2$ are false is the only case in four posible cases (given that there are two possible truth values for $q_1$ and two possible values for $q_2$) in which the proposition $q_1 \lor q_2$ is considered false. In the other three possible cases (in which at least one of the two propositions $q_1$ and $q_2$ is true) $q_1 \lor q_2$ is considered true.

In the first row of the membership table in figure \ref{f2b} the numerical sequence shown is 0, 0, 0. It should be interpreted as follows: If any element whatsoever belonging to the universal set $\mathbb U$ considered belongs neither to $C_1$ nor to $C_2$, then that element does not belong to the set union of $C_1$ and $C_2$ (that is, to $C_1 \cup C_2$) either. Note that the case in which any element whatsoever of $\mathbb U$ belongs neither to $C_1$ nor to $C_2$ is the only case of the four possible ones (given that there are two possible membership values $(0,1)$ for that element to belong to $C_1$ and two possible membership values (0, 1) for that same element to belong to $C_2$), in which that element does not belong to $C_1 \cup C_2$. In the other three possible cases (in which that element belongs to at least one of the sets $C_1$ and $C_2$), that element does belong to $C_1 \cup C_2$.

From the viewpoint of the presence of zeros and ones, the truth table for $q_1 \lor q_2$ is equal to the membership table for $C_1 \cup C_2$. This result makes it possible to establish the following correspondences: 1) between $q_1$ and $C_1$, 2) between $q_2$ and $C_2$, 3) between the operator of disjunction in propositional calculus and the operator of union in set theory, and 4) between the proposition $q_1 \lor q_2$ and the set $C_1 \cup C_2$.

In propositional calculus there is also another type of operation of disjunction: that denominated ``exclusive disjunction''. When dealing with a topic in which both types of disjunction are used, the disjunction considered above is usually denominated `inclusive disjunction'', to prevent confusion.

The operator of \textit{exclusive disjunction} of two propositions $q_1$ and $q_2$ is symbolized as $\dot{\lor}$. The proposition $q_1 \dot{\lor} q_2$ is read as ``$q_1$ exclusive or $q_2$''.

The operator of \textit{exclusive union} of two sets $C_1$ and $C_2$ is symbolized as $\dot{\cup}$. The set $C_1 \dot{\cup} C_2$ is the set resulting from the exclusive union of the sets $C_1$ and $C_2$.

The truth table for the proposition $q_1 \dot{\lor} q_2$ is shown in figure \ref{f3a}, and the membership table for the set $C_1 \dot{\cup} C_2$ is shown in figure \ref{f3b}.

\begin{figure}[H]
\centering
\subfloat[Truth table for $q_1 \dot{\lor} q_2$]{
\hspace{1.6in}
\begin{tabular}{c|c||c|}
$q_1$ &  $q_2$ & $q_1 \dot{\lor} q_2$ \\
\midrule
0 & 0 & 0   \\  
0 & 1 & 1   \\  
1 & 0 & 1   \\  
1 & 1 & 0   \\  
\end{tabular}
\hspace{2in}
\label{f3a}
}
\\
\subfloat[Membership table for the set $C_1 \dot{\cup} C_2$]{
\hspace{1.5in}
\begin{tabular}{c|c||c|}
$C_1$ &  $C_2$ & $C_1 \dot{\cup} C_2$ \\
\midrule
0 & 0 & 0   \\  
0 & 1 & 1   \\  
1 & 0 & 1   \\  
1 & 1 & 0   \\  
\end{tabular}
\hspace{2in}
\label{f3b}}%
\caption{a) Truth table for $q_1 \dot{\lor} q_2$, and b) Membership table for $C_1 \dot{\cup} C_2$}
\label{f3}
\end{figure}

In the truth table in figure \ref{f3a} it can be seen that $q_1 \dot{\lor} q_2$ is true if $q_1$ and $q_2$ have different truth values. This occurs in the case in which $q_1$ is false and $q_2$ is true (corresponding to the second row of the truth table) and in the case in which $q_1$ is true and $q_2$ is false (corresponding to the third row of the truth table). Both in the case in which $q_1$ and $q_2$ are false (shown in the first row of the truth table), and in the case in which $q_1$ and $q_2$ are true (shown in the fourth row of the truth table), $q_1 \dot{\lor} q_2$ is false.

In the membership table in figure \ref{f3b} it can be seen that only in the cases in which any element whatsoever of the universal set $\mathbb U$ belongs to one of the two sets $C_1$ and $C_2$, and does not belong to the other (those corresponding to the second and third rows of that membership table), does that element belong to $C_1 \dot{\cup} C_2$. Both if any element whatsoever of $\mathbb U$ belongs neither to $C_1$ nor to $C_2$ (as in the case shown in the first row of that membership table) and if any element whatsoever of $\mathbb U$ belongs to $C_1$ and to $C_2$ (as in the case shown in the fourth row of the membership table), that element does not belong to $C_1 \dot{\cup} C_2$.

From the viewpoint of the presence of zeros and ones, the truth table for $q_1 \dot{\lor} q_2$ is equal to the membership table for $C_1 \dot{\cup} C_2$. This result makes it possible to establish the following correspondences: 1) between $q_1$ and $C_1$, 2) between $q_2$ and $C_2$, 3) between the operator of exclusive disjunction in propositional calculus and the operator of exclusive union in set theory, and 4) between the proposition $q_1 \dot{\lor} q_2$ and the set $C_1 \dot{\cup} C_2$.

The operator of \textit{conjunction} of two propositions $q_1$ and $q_2$ is symbolized as $\land$. The proposition $q_1 \land q_2$ is read as ``$q_1$ and $q_2$''.

The operator of \textit{intersection} of two sets $C_1$ and $C_2$ is symbolized as $\cap$. The set $C_1 \cap C_2$ is the set resulting from the intersection of the sets $C_1$ and $C_2$.

The truth table for the proposition $q_1 \land q_2$ is shown in figure \ref{f4a}, and the membership table for the set $C_1 \cap C_2$ is shown in figure \ref{f4b}.

\begin{figure}[H]
\centering
\subfloat[Truth table for $q_1 \land q_2$]{
\hspace{1.6in}
\begin{tabular}{c|c||c|}
$q_1$ &  $q_2$ & $q_1 \land q_2$  \\
\midrule
0 & 0 & 0   \\  
0 & 1 & 0   \\  
1 & 0 & 0   \\  
1 & 1 & 1   \\  
\end{tabular}
\hspace{2in}
\label{f4a}
}
\\
\subfloat[Membership table for $C_1 \cap C_2$]{
\hspace{1.5in}
\begin{tabular}{c|c||c|}
$C_1$ &  $C_2$ & $C_1 \cap C_2$  \\
\midrule
0 & 0 & 0   \\  
0 & 1 & 0   \\  
1 & 0 & 0   \\  
1 & 1 & 1   \\  
\end{tabular}
\hspace{2in}
\label{f4b}}%
\caption{a) Truth table for $q_1 \land q_2$; and b) Membership table for $C_1 \cap C_2$}
\label{f4}
\end{figure}

Note that in the truth table shown in figure \ref{f4a} the fourth row -- with the numerical sequence 1, 1, 1 -- corresponds to the only case (that in which both $q_1$ and $q_2$ are true) in which $q_1 \land q_2$ is true. In the case in which $q_1$ and $q_2$ are false (corresponding to the first row of that truth table), as in the cases in which only one of the propositions $q_1$ and $q_2$ is true (corresponding to the second and third rows of that table), $q_1 \land q_2$ is false.

Observe that in the membership table shown in figure \ref{f4b} the fourth row -- with the numerical sequence 1, 1, 1 -- corresponds to the only case in which any element whatsoever of $\mathbb U$ belongs both to $C_1$ and $C_2$ (So that element belongs to $C_1 \cap C_2$.). In the case in which any element whatsoever of $\mathbb U$ belongs neither to $C_1$ nor to $C_2$ (corresponding to the first row of the membership table), as in the cases in which any element whatsoever of $\mathbb U$ belongs to only one of the sets $C_1$ and $C_2$ (corresponding to the second and third rows of that membership table), that element does not belong to $C_1 \cap C_2$.

From the viewpoint of the presence of zeros and ones, the truth table for $q_1 \land q_2$ is equal to the membership table for $C_1 \cap C_2$. This result makes it possible to establish the following correspondences: 1) between $q_1$ and $C_1$, 2) between $q_2$ and $C_2$, 3) between the operator of conjunction in propositional calculus and the operator of intersection in set theory, and 4) between the proposition $q_1 \land q_2$ and the set $C_1 \cap C_2$.

The operator of \textit{material implication} in propositional calculus will be symbolized as $\to$. Emphasis must be given to an important difference between material implication and other operations in propositional calculus, such as inclusive disjunction, exclusive disjunction, and conjunction; that is, the propositions $q_1 \lor q_2$, $q_1 \dot{\lor} q_2$, and $q_1 \land q_2$ have the same meaning and the same truth tables, respectively, as the propositions $q_2 \lor q_1$, $q_2 \dot{\lor} q_1$, and $q_2 \land q_1$. However, $q_1 \to q_2$, which is read as ``if $q_1$, then $q_2$'' (or ``$q_1$ materially implies $q_2$''), has a meaning and a truth table different from those of the proposition $q_2 \to q_1$, which is read as ``if $q_2$, then $q_1$'' (or as ``$q_2$ materially implies $q_1$''). Both $q_1 \to q_2$ and $q_2 \to q_1$ are considered propositions of a conditional character.

The proposition $q_1$ is denominated the antecedent of the proposition $q_1 \to q_2$, and $q_2$ is denominated the consequent of that proposition. The proposition $q_2$ is denominated the antecedent of the proposition $q_2 \to q_1$, and $q_1$ is denominated the consequent of that proposition.

The operator of material implication in set theory is symbolized as $\naturalto$. The membership table for the set $C_1 \naturalto C_2$ is different from that of $C_2 \naturalto C_1$. The consideration of these membership tables will make it possible to characterize the sets $C_1 \naturalto C_2$ and $C_2 \naturalto C_1$ unambiguously.

The truth tables of $q_1 \to q_2$, and $q_2 \to q_1$ are shown, respectively in figures \ref{f5a} and \ref{f5b}. The membership tables of $C_1 \naturalto C_2$ and $C_2 \naturalto C_1$ are shown, respectively in figures \ref{f5c} and \ref{f5d}. 
\begin{figure}[H]
\centering
\subfloat[Truth table for $q_1 \to q_2$]{
\begin{tabular}{c|c||c|}
$q_1$ &  $q_2$ & $q_1 \to q_2$  \\
\midrule
0 & 0 & 1 \\  
0 & 1 & 1  \\  
1 & 0 & 0  \\  
1 & 1 & 1 \\  
\end{tabular}
\hspace{0.5in}
\label{f5a}
}
\subfloat[Truth table for $q_2 \to q_1$]{
\begin{tabular}{c|c||c|}
$q_1$ &  $q_2$  & $q_2 \to q_1$ \\
\midrule
 0 & 0 & 1 \\  
 0 & 1 & 0 \\  
 1 & 0 & 1 \\  
 1 & 1 & 1 \\  
\end{tabular}
\hspace{1in}
\label{f5b}
}
\\
\subfloat[Membership table for $C_1 \naturalto C_2$]{
\begin{tabular}{c|c||c|} 
$C_1$ &  $C_2$ & $C_1 \naturalto C_2$ \\
\midrule
0 & 0 & 1   \\  
0 & 1 & 1   \\  
1 & 0 & 0  \\  
1 & 1 & 1   \\    
\end{tabular}
\hspace{0.5in}
\label{f5c}
}
\subfloat[Membership table for $C_2 \naturalto C_1$]{
\begin{tabular}{c|c||c|} 
 $C_1$ &  $C_2$ & $C_2 \naturalto C_1$ \\
\midrule
  0 & 0 & 1 \\  
  0 & 1 & 0 \\  
 1 & 0 & 1 \\  
 1 & 1 & 1 \\    
\end{tabular}
\hspace{1in}
\label{f5d}
}
\caption{a) Truth table for $q_1 \to q_2$, and b) Truth table for $q_2 \to q_1$;\\ c) Membership table for $C_1 \naturalto C_2$, and d) Membership table for $C_2 \naturalto C_1$}
\label{f5}
\end{figure}

Note in the truth table in figure \ref{f5a} that the only case in which the proposition $q_1 \to q_2$ is considered false is that corresponding to the third row of the table, with the numerical sequence 1, 0, 0. If $q_1$ is true and $q_2$ is false, then $q_1 \to q_2$ is false. In the other three possible cases considered in this truth table, $q_1 \to q_2$ is considered true.

Note in the truth table in figure \ref{f5b} that the only case in which the proposition $q_2 \to q_1$ is false is that of the second row of that table, with the numerical sequence 0, 1, 0: If $q_2$ is true and $q_1$ is false, then $q_2 \to q_1$ is false. In the other three possible cases considered in this table, $q_2 \to q_1$ is true.

Observe in the membership table in figure \ref{f5c} that the only case in which any element whatsoever of the set $\mathbb U$ considered does not belong to the set $C_1 \naturalto C_2$ is that of the third row of the table, with the numerical sequence 1, 0, 0: If that element belongs to $C_1$ and does not belong to $C_2$, then that element does not belong to $C_1 \naturalto C_2$. In the other three possible cases shown in that membership table, that element does belong to $C_1 \naturalto C_2$.

Note in the membership table in figure \ref{f5d} that the only case in which any element whatsoever of the set $\mathbb U$ considered does not belong to the set $C_2 \naturalto C_1$ is that of the second row of the table, with the numerical sequence 0, 1, 0: If that element belongs to $C_2$ and does not belong to $C_1$, then that element does not belong to $C_2 \naturalto C_1$. In the other three possible cases considered in that membership table, that element does belong to $C_2 \naturalto C_1$.

Therefore, the following correspondences may be established: 1) between $q_1$ and $C_1$, as well as between $q_2$ and $C_2$; 2) between the operator of material implication in propositional calculus and the operator of material implication in set theory; and 3) between $q_1 \to q_2$ and $C_1 \naturalto C_2$, as well as between $q_2 \to q_1$ and $C_2 \naturalto C_1$.

The operator of \textit{material bi-implication}, or of logical equivalence, in propositional calculus will be symbolized as $\longleftrightarrow$. The operator of material bi-implication in set theory will be symbolized as $\naturaltolr$.

The truth table for the proposition $q_1 \longleftrightarrow q_2$ is shown in figure \ref{f6a}. The membership table for the set $C_1 \naturaltolr C_2$ is shown in figure \ref{f6b}.

\begin{figure}[H]
\centering
\subfloat[Truth table for the proposition $q_1 \longleftrightarrow q_2$]{
\hspace{1.5in}
\begin{tabular}{c|c||c|}
$q_1$ &  $q_2$ & $q_1 \longleftrightarrow q_2$  \\
\midrule
0 & 0 & 1  \\  
0 & 1 & 0   \\  
1 & 0 & 0   \\  
1 & 1 & 1   \\  
\end{tabular}
\hspace{2in}
\label{f6a}
}
\\
\subfloat[Membership table for the set $C_1 \naturaltolr C_2$]{
\hspace{1.4in}
\begin{tabular}{c|c||c|}
$C_1$ &  $C_2$ & $C_1 \naturaltolr C_2$  \\
\midrule
0 & 0 & 1  \\  
0 & 1 & 0   \\  
1 & 0 & 0   \\  
1 & 1 & 1   \\  
\end{tabular}
\hspace{2in}
\label{f6b}}%
\caption{a) Truth table for $q_1 \longleftrightarrow q_2$; and b) Membership table for the set $C_1 \naturaltolr C_2$}
\label{f6}
\end{figure}

Note in the truth table in figure \ref{f6a} that there are two cases in which the proposition $q_1 \longleftrightarrow q_2$ is true: those in which $q_1$ and $q_2$ have the same truth value. The case in which $q_1$ and $q_2$ are false corresponds to the first row of that truth table, with the numerical sequence 0, 0, 1. The case in which $q_1$ and $q_2$ are true corresponds to the fourth row of that truth table, with the numerical sequence 1, 1, 1. In the other two possible cases considered in that truth table, in which $q_1$ and $q_2$ have different truth values, the proposition $q_1 \longleftrightarrow q_2$ is false.

Observe in the membership table in figure \ref{f6b} that there are two cases in which any element whatsoever of the universal set $\mathbb U$ considered belongs to the set $C_1 \naturaltolr C_2$: the case in which that element belongs neither to $C_1$ nor to $C_2$ and the case in which that element belongs both to $C_1$ and to $C_2$. The case in which that element belongs neither to $C_1$ nor to $C_2$ corresponds to the first row of that membership table, with the numerical sequence 0, 0, 1. The case in which any element whatsoever of $\mathbb U$ belongs both to $C_1$ and to $C_2$ corresponds to the fourth row of that membership table, with the numerical sequence 1, 1, 1. In the other two possible cases considered in that membership table, in which any element whatsoever of the universal set $\mathbb U$ belongs to one of the two sets and not to the other, that element does not belong to the set $C_1 \naturaltolr C_2$.

From the viewpoint of the presence of zeros and ones, the truth table in figure \ref{f6a} is equal to the membership table in figure \ref{f6b}. This equality makes it possible to establish the following correspondences: 1) between $q_1$ and $C_1$, 2) between $q_2$ and $C_2$, 3) between the operator of material bi-implication in propositional calculus and the operator of material bi-implication in set theory, and 4) between the proposition $q_1 \longleftrightarrow q_2$ and the set $C_1 \naturaltolr C_2$ .

The \textit{Sheffer stroke} operator, or ``nand,'' in propositional calculus will be symbolized as $\uparrow$. The \textit{Sheffer stroke} operator, or ``nand,'' in set theory will be symbolized as $\uparrowbarred$.

The truth table for the proposition $q_1 \uparrow q_2$ is shown in figure \ref{f7a}. The membership table for the set $C_1 \uparrowbarred C_2$ is shown in figure \ref{f7b}.

\begin{figure}[H]
\centering
\subfloat[Truth table for the proposition $q_1 \uparrow q_2$]{
\hspace{1.6in}
\begin{tabular}{c|c||c|}
$q_1$ &  $q_2$ & $q_1 \uparrow q_2$ \\
\midrule
0 & 0 & 1   \\  
0 & 1 & 1   \\  
1 & 0 & 1   \\  
1 & 1 & 0   \\  
\end{tabular}
\hspace{2in}
\label{f7a}
}
\\
\subfloat[Membership table for the set $C_1 \uparrowbarred C_2$]{
\hspace{1.5in}
\begin{tabular}{c|c||c|}
$C_1$ &  $C_2$ & $C_1 \uparrowbarred C_2$  \\
\midrule
0 & 0 & 1   \\  
0 & 1 & 1   \\  
1 & 0 & 1   \\  
1 & 1 & 0   \\  
\end{tabular}
\hspace{2in}
\label{f7b}}%
\caption{a) Truth table for the proposition $q_1 \uparrow q_2$; and b) Membership table for the set $C_1 \uparrowbarred C_2$}
\label{f7}
\end{figure}

Note that the truth table in figure \ref{f7a}, the only case in which the proposition $q_1 \uparrow q_2$ is false is that in which both $q_1$ and $q_2$ are true. This case corresponds to the fourth row of that truth table, with the numerical sequence 1, 1, 0. In the case in which both $q_1$ and $q_2$ are false, corresponding to the first row of that truth table, with the numerical sequence 0, 0, 1, and in the two cases in which $q_1$ and $q_2$ have different truth values, corresponding to the second and third rows of that truth table, the proposition $q_1 \uparrow q_2$ is true.

Observe that in the membership table in figure \ref{f7b} that the only case in which any element whatsoever of the universal set $\mathbb U$ considered does not belong to the set $C_1 \uparrowbarred C_2$ is that in which that element belongs both to $C_1$ and to $C_2$. In the case in which any element whatsoever of $\mathbb U$ belongs neither to $C_1$ nor to $C_2$ (corresponding to the first row in that membership table, with the numerical sequence 0, 0, 1), and in the two possible cases in which any element whatsoever of $\mathbb U$ belongs to only one of the sets $C_1$ and $C_2$ (corresponding to the second and third rows of that membership table) that element belongs to $C_1 \uparrowbarred C_2$.

From the viewpoint of the presence of zeros and ones, the truth table in figure \ref{f7a} is equal to the membership table in figure \ref{f7b}. This equality makes it possible to establish the following correspondences: 1) between $q_1$ and $C_1$, 2) between $q_2$ and $C_2$, 3) between the Sheffer stroke operator in propositional calculus and the Sheffer stroke operator in set theory, and 4) between the proposition $q_1 \uparrow q_2$ and the set $C_1 \uparrowbarred C_2$.

The \textit{Peirce's arrow} operator -- or ``nor'' -- in propositional calculus will be symbolized as $\downarrow$. The \textit{Peirce's arrow} operator -- or ``nor''-- in set theory will be symbolized as $\downarrowbarred .$

The truth table for the proposition $q_1 \downarrow q_2$ is shown in figure \ref{f8a}. The membership table for the set $C_1 \downarrowbarred C_2$ is shown in figure \ref{f8b}.

\begin{figure}[H]
\centering
\subfloat[Truth table for the proposition $q_1 \downarrow q_2$]{
\hspace{1.6in}
\begin{tabular}{c|c||c|}
$q_1$ &  $q_2$ & $q_1 \downarrow q_2$ \\
\midrule
0 & 0 & 1   \\  
0 & 1 & 0   \\  
1 & 0 & 0   \\  
1 & 1 & 0   \\  
\end{tabular}
\hspace{2in}
\label{f8a}
}
\\
\subfloat[Membership table for the set $C_1 \downarrowbarred C_2$]{
\hspace{1.5in}
\begin{tabular}{c|c||c|}
$C_1$ &  $C_2$ & $C_1 \downarrowbarred C_2$  \\
\midrule
0 & 0 & 1   \\  
0 & 1 & 0   \\  
1 & 0 & 0   \\  
1 & 1 & 0   \\  
\end{tabular}
\hspace{2in}
\label{f8b}}%
\caption{a) Truth table for the proposition $q_1 \downarrow q_2$; and b) Membership table for the set $C_1 \downarrowbarred C_2$}
\label{f8}
\end{figure}

Note that in the truth table in \ref{f8a}, the only case in which the proposition $q_1 \downarrow q_2$ is true is that in which both $q_1$ and $q_2$ are false. This case corresponds to the first row of that truth table, with the numerical sequence 0, 0, 1. In the case in which both $q_1$ and $q_2$ are true (corresponding to the fourth row of the truth table, with the numerical sequence 1, 1, 0), and in the two cases in which $q_1$ and $q_2$ have different truth values (corresponding to the second and third rows of that truth table), the proposition $q_1 \downarrow q_2$ is false.

Observe in the membership table in figure \ref{f8b} that the only case in which any element whatsoever of the universal set $\mathbb U$ considered belongs to the set $C_1 \downarrowbarred C_2$ is that in which that element belongs neither to $C_1$ nor to $C_2$. In the case in which any element whatsoever of $\mathbb U$ belongs both to $C_1$ and to $C_2$ (corresponding to the fourth row in that membership table, with the numerical sequence 1, 1, 0) and in the two possible cases in which any element whatsoever of $\mathbb U$ belongs to only one of the sets $C_1$ and $C_2$ (corresponding to the second and third rows of that membership table) that element does not belong to the set $C_1 \downarrowbarred C_2$.

From the viewpoint of the presence of zeros and ones, the truth table in figure \ref{f8a} is equal to the membership table in figure \ref{f8b}. This equality makes it possible to establish the following correspondences: 1) between $q_1$ and $C_1$, 2) between $q_2$ and $C_2$, 3) between the Peirce's arrow operator in propositional calculus and the Peirce's arrow operator in set theory, and 4) between the proposition $q_1 \downarrow q_2$ and the set $C_1 \downarrowbarred C_2$.

From the standpoint of the information required to determine which of the elements belonging to a set obtained as the result of any operation in set theory, two types of operations may be distinguished in this theory.

Suppose that the operation of union of two sets $C_1$ and $C_2$ is carried out. To obtain the set resulting from that operation $C_1 \cup C_2$, it is not indispensable to entirely know the universal set $\mathbb U$ within whose frame $C_1$ and $C_2$ were characterized. Indeed, to compute $C_1 \cup C_2$ it suffices to know which elements of that $\mathbb U$ belong to $C_1$ and which elements of that $\mathbb U$ belong to $C_2$. As seen above, any element belonging to at least one of the two sets $C_1$ and $C_2$ belongs to $C_1 \cup C_2$. Other elements belonging to $\mathbb U$ that belong neither to $C_1$ nor to $C_2$ whose knowledge is not indispensable to compute $C_1 \cup C_2$ can exist.

Another operation of the same type is that of the intersection of two sets $C_1$ and $C_2$. Indeed, the elements of 
 $\mathbb U$ (that is, those belonging to the $\mathbb U$ considered) which belong both to $C_1$ and $C_2$, belong to the result of the operation of the intersection of $C_1$ and $C_2$ (that is, to the set $C_1 \cap C_2$). Other elements belonging to $\mathbb U$ that belong neither to $C_1$ nor to $C_2$ whose knowledge is not indispensable to compute $C_1 \cap C_2$ can exist.
 
The operation of exclusive union is also an operation of the same type in set theory. Indeed, to compute $C_1 \dot{\cup} C_2$ (that is, to determine which elements belong to that set), it suffices to know which elements belong to 
$C_1$ and which elements belong to $C_2$. If an element belonging to one of these two sets does not belong to the other, then that element belongs to $C_1 \dot{\cup} C_2$.

However, if consideration is given, for example, to the operation of complementation in set theory, it does not suffice to know which elements belong to any set $C$ whatsoever to be able to compute its complement $\xvec{C}$. Indeed, to determine which elements belong to $\xvec{C}$ (that is, all elements of $\mathbb U$ which do not belong to $C$), it is necessary to know which elements belong to $\mathbb U$. The operation of the complementation of a set, in the sense specified above, is of a type different from those already considered.

Another operation of this second type (that is, such that to be able to carry it out, it is necessary to know which $\mathbb U$ is considered, and not only the sets on which it is acting), is that of material implication. In effect, if any element whatsoever of $\mathbb U$ belongs neither to a set $C_1$ nor to a set $C_2$, both characterized within the frame of some $\mathbb U$, then that element does belong to $C_1 \naturalto C_2$ as well as to $C_2 \naturalto C_1$. It is clear that the mere knowledge of $C_1$ and $C_2$ (that is, of which elements belong to each of these two sets) does not suffice to know which elements of $\mathbb U$ belong neither to $C_1$ nor to $C_2$.

Other operations in set theory, considered above, also correspond to the latter type. Indeed, in each membership table corresponding to each of the sets $C_1 \naturaltolr C_2$, $C_1 \uparrowbarred C_2$ and $C_1 \downarrowbarred C_2$, there is a row with the numerical sequence 0, 0, 1. That is, if any element whatsoever of the $\mathbb U$ considered belongs neither to $C_1$ nor to $C_2$, that element does belong to one of the sets $C_1 \naturaltolr C_2$, $C_1 \uparrowbarred C_2$ and $C_1 \downarrowbarred C_2$. To be able to compute each of these last three sets, it does not suffice, therefore, to know which elements belong to $C_1$ and which belong to $C_2$; it is also necessary to know what that $\mathbb U$ is (that is, to know all the elements belonging to the $\mathbb U$ considered).

\section{How Parentheses Will Be Used for Propositional Calculus and Set Theory}

This section concerns one aspect of the notation used in this article.

Suppose that one wants to obtain an exclusive disjunction of the proposition $q_1 \lor q_2$ (that is, the inclusive disjunction of $q_1$ and $q_2$) and the proposition $q_1 \to q_3$ (that is, the material implication of which the proposition $q_1$ is the antecedent and the proposition $q_3$ is the consequent). To eliminate any ambiguity regarding the exclusive disjunction referred to, the specific operation can be expressed as \noindent $(q_1 \lor q_2)\dot{\lor}(q_1 \to q_3)$.

It can be seen that the proposition $q_1 \lor q_2$ was specified in parentheses, as was the proposition $q_1 \to q_3$.

The conjunction of the exclusive disjunction above $(q_1 \lor q_2)\dot{\lor}(q_1 \to q_3)$ and the proposition $(q_3 \to q_4)$ can be specified unambiguously as follows:
\begin{equation}
((q_1 \lor q_2)\dot{\lor}(q_1 \to q_3))\land(q_3 \to q_4) 
\end{equation}

Note that a set corresponding to the proposition $(q_1 \lor q_2)\dot{\lor}(q_1 \to q_3)$ is the exclusive union of the sets 
$C_1 \cup C_2$ and $C_1 \to C_3$. To prevent any ambiguity, that union may be expressed as follows:
\begin{equation}
(C_1 \cup C_2)\dot{\cup}(C_1 \naturalto C_3) 
\end{equation}

A set corresponding to the material implication of $q_3 \to q_4$ is the following:
\begin{equation}
C_3 \naturalto C_4 
\end{equation}

Set (4) corresponding to proposition (1) is the intersection of the sets specified in (2) and (3):
\begin{equation}
((C_1 \cup C_2)\dot{\cup}(C_1 \naturalto C_3))\cap(C_3 \naturalto C_4)  
\end{equation}

To begin to characterize the use to be made of parentheses in this article in propositional calculus, suppose first that an operation is carried out with two propositions and another proposition is obtained. If either of these two propositions is a function of two or more propositions (that is, in its expression there are two or more propositions), then it will be expressed in parentheses. Thus, for example, consider once more the following two propositions: a) $q_1 \lor q_2$ and b) $q_1 \to q_3$. Note that $q_1 \lor q_2$ is a function of the propositions $q_1$ and $q_2$ and that $q_1 \to q_3$ is a function of the propositions $q_1$ and $q_3$. Therefore, if an operation -- such as the operation of exclusive disjunction -- is carried out with the propositions $q_1 \lor q_2$ and $q_1 \to q_3$, the proposition resulting from that operation will be expressed as \noindent ($q_1 \lor q_2)\dot{\lor}(q_1 \to q_3$).

If, with the latter proposition (which is a function of the three propositions $q_1$, $q_2$ and $q_3$), another operation -- such as the operation of conjunction -- is carried out with the proposition $q_3 \to q_4$, which is a function of the two propositions $q_3$ and $q_4$, the proposition resulting from that operation of conjunction will be expressed as follows:
\begin{equation*}
((q_1 \lor q_2)\dot{\lor}(q_1 \to q_3)) \land (q_3 \to q_4)  
\tag{1}
\end{equation*}

As mentioned above, the set corresponding to proposition (1) is the following:
\begin{equation*}
((C_1 \cup C_2)\dot{\cup}(C_1 \naturalto C_3))\cap(C_3 \naturalto C_4)  
\tag{4}
\end{equation*}

Proposition (1) and set (4) are isomorphic: $q_1$, $q_2$, $q_3$ and $q_4$ correspond respectively to $C_1$, $C_2$, $C_3$ and $C_4$; and the operators of inclusive disjunction ($\lor$), exclusive disjunction ($\dot{\lor}$), material implication, ($\to$) and conjunction ($\land$) in propositional calculus correspond, respectively, to the operators of union ($\cup$), exclusive union ($\dot{\cup}$), material implication ($\naturalto$) and intersection ($\cap$) in set theory.

To complete the specification of the use of parentheses in this article in propositional calculus and in set theory, it is specified that propositions which are functions of a sole proposition but in which its negation is also present will be expressed in parentheses when carrying out operations. Thus, for example, suppose that one desires to express the exclusive disjunction of $q_1 \lor \overline{q}_1$ and $q_2 \to q_3$. This exclusive disjunction will be expressed as follows:
\begin{equation}
(q_1 \lor \overline{q}_1) \dot{\lor}(q_2 \to q_3) 
\end{equation}

Observe that $q_1 \lor \overline{q}_1)$ is true, regardless of the truth value of $q_1$.

The set corresponding to proposition (5) is the following: $(C_1 \cup \xvec{C}_1) \dot{\cup}(C_2 \naturalto C_3)$.

Note that $C_1 \cup \xvec{C}_1$ is equal to the universal set $\mathbb U$ within whose frame $C_1$ has been characterized: $C_1 \cup \xvec{C}_1 = \mathbb U$.

Another example in which parentheses are used for a proposition which is the function of a sole proposition is the following: Admit that $q_1$ materially implies the conjunction $q_2$ and $\overline{q}_2$:
\begin{equation}
q_1 \to (q_2 \land \overline{q}_2) 
\end{equation}

In the material implication above its consequent has been shown in parentheses.

Note that $q_2 \land \overline{q}_2$ is false, regardless of the truth value of $q_2$.

The set corresponding to the proposition (6) is $C_1 \naturalto (C_2 \cap \xvec{C}_2)$.

Note that $C_2 \cap \xvec{C}_2$ is equal to the empty set: $C_2 \cap \xvec{C}_2 = \varnothing$.

\section{Number of Rows in the Truth Table of a Proposition that Is a Function of $n$ Propositions, for $n = 1, 2, 3,\dots$, and Number of Logical Functions of that Proposition}

Recall that the number of rows in a truth table of a proposition that is a function of $n$ propositions, for $n = 1, 2, 3,\dots$, is $2^n$. In effect, given that each of those $n$ propositions can have two truth values (true or false), there are $2^n$ cases possible for assigning truth values to those $n$ propositions, such that each case differs from the rest.

The characterization of the function considered of those $n$ propositions is carried out by establishing in the column corresponding to that function, in that truth table, for each one of the $2^n$ cases, whether the function considered -- also a proposition -- is true or false. It is then observed that $2^{(2^{n})}$ different possible truth values can be assigned in the column corresponding to the function. In other words, there are $2^{(2^{n})}$ logical functions of each proposition, which is a function of $n$ propositions.

Thus, for example, if a proposition is a function of two propositions, the corresponding truth table has 4 (that is, $2^2$) rows and there are 16 -- that is, $2^{(2^{2})}$ possible logical functions of that proposition.

Consider a proposition that is a function of three propositions. In this case, the corresponding truth table has 8 (that is, $2^3$) rows and there are 256 (that is, $2^{(2^{3})}$) possible logical functions of that proposition.

\section{Certain Laws -- or Tautologies -- in Propositional Calculus and the Corresponding Laws in Set Theory}

If a proposition which is a function of $n$ propositions -- $n=1, 2, 3, \dots$ -- is true regardless of the truth values of each of those $n$ propositions, then it is considered a law, or tautology, in propositional calculus. For each law in propositional calculus there is set that is isomorphic to it which is equal to the universal set within whose frame that set has been characterized.

In this section attention is given to several laws of propositional calculus and the corresponding sets which are expressions of the universal sets within whose frames those sets have been characterized.

In figure \ref{f9a} consideration is given to the law of propositional calculus $q \lor \overline{q}$. This law is denominated ``law of the excluded middle''. In the truth table presented in this figure it can be observed that both if $q$ if true and if $q$ is false, the proposition $q \lor \overline{q}$ is true.

Given any universal set $\mathbb U$ and any set $C$ characterized within that frame, the set $C \cup \xvec{C}$ resulting from the operation of union of $C$ and its complement $ \xvec{C}$ is a) isomorphic to the proposition 
$q \lor \overline{q}$, and b) equal to the universal set $C \cup \xvec{C} = \mathbb U$. The membership table for the set $C \cup \xvec{C}$ is presented in figure \ref{f9b}.

\begin{figure}[H]
\centering
\subfloat[Truth table for the proposition $q_1 \lor \overline q$.]{
\hspace{1.7in}
\begin{tabular}{c|c||c|}
$q$ & $\overline{q}$  &  $q\lor \overline{q}$\\
\midrule
0 & 1 & 1  \\  
1 & 0  & 1 \\  
\end{tabular}
\hspace{2in}
\label{f9a}
}
\\
\subfloat[Membership table for the set $C \cup \xvec{C}$.]{
\hspace{1.6in}
\begin{tabular}{c|c||c|}
$C$ & ${\xvec{C}}$  & $C \cup \xvec{C}$ \\
\midrule
0 & 1   & 1 \\  
1 & 0   & 1 \\  
\end{tabular}
\hspace{2in}
\label{f9b}}
\caption{a) Truth table for the proposition $q\lor \overline{q} $, and b) Membership table for the set $C \cup \xvec{C}$}
\label{f9}
\end{figure}

As observed in figure \ref{f9a}, regardless of the truth value of the proposition $q$ -- considered false in the first row of the corresponding truth table and true in the second row -- $q\lor \overline{q}$ is true. Precisely for this reason, $q\lor \overline{q}$ is considered a law in propositional calculus. Likewise, as observed in  figure \ref{f9b}, regardless of whether any element of $\mathbb U$ within whose frame $C$ was characterized does not belong to $C$ (as seen in the first row of the corresponding membership table), or does belong to $C$ (as seen in the second row of the table), that element belongs to $C \cup \xvec{C}$. Precisely for this reason, it is admitted that $C \cup \xvec{C}$ is equal to the universal set considered.

In this section it can be seen that in the column corresponding to each proposition whose truth table is presented there are only numerical values equal to one (1). In other words, regardless of the truth values of the different propositions $q_1, q_2, q_3, \dots$ of which the proposition whose truth table is considered is a function, this proposition is true. Precisely for this reason the proposition whose truth table is presented is considered a law, or tautology, of propositional calculus. Likewise, it can be observed in the column of the membership table corresponding to each set that is isomorphic to a law in propositional calculus there are only numerical values equal to one (1). In other words, regardless of whether each element of the universal set considered belongs or does not belong to each one of the sets $C_1, C_2, C_3, \dots$ of which the set whose membership table is presented is a function, that element belongs to this set.
Precisely for this reason this set is equal to the universal set.

The law of propositional calculus $((q_1 \to q_2) \land q_1) \to q_2$ is known as ``modus ponendo ponens.'' It can be considered as a way of reasoning such that if there is a conditional relation $q_1 \to q_2$ and if it is affirmed (\textit{ponendo}) that the antecedent $q_1$ is true, it is inferred and affirmed (\textit{ponens}) that the consequent $q_2$ is true.

The truth table of the law of propositional calculus $((q_1 \to q_2) \land q_1) \to q_2$ is presented in figure \ref{f10a}, and the membership table of the set that is isomorphic to that law $((C_1 \naturalto C_2) \cap C_1) \naturalto C_2$ is presented in \ref{f10b}.

\begin{figure}[H]
\centering
\subfloat[	Truth table for the proposition $((q_1 \to q_2) \land q_1) \to q_2 $]{
\begin{tabular}{c|c|c|c||c|} 
$q_1$ &  $q_2$ & $q_1 \to q_2$  & $ (q_1 \to q_2) \land q_1 $ & $((q_1 \to q_2) \land q_1) \to q_2$\\
\midrule
0 & 0 & 1 & 0 & 1  \\ 
0 & 1 & 1 & 0 & 1   \\  
1 & 0 & 0 & 0 & 1   \\  
1 & 1 & 1 & 1 & 1   \\  
\end{tabular}
\label{f10a}
}
\\
\subfloat[Membership table for the set $((C_1 \naturalto C_2) \cap C_2 ) \naturalto C_1$; $(((C_1 \naturalto C_2) \cap C_1) \naturalto C_2) = \mathbb U$]{
\hspace{-.2in}
\begin{tabular}{c|c|c|c||c|}
$C_1$ &  $C_2$ & $C_1 \naturalto C_2$ & $(C_1 \naturalto C_2) \cap C_1$  & $((C_1 \naturalto C_2) \cap C_1) \naturalto C_2$\\
\midrule
0 & 0 & 1 & 0 & 1  \\  
0 & 1 & 1 & 0 & 1   \\  
1 & 0 & 0 & 0 & 1   \\  
1 & 1 & 1 & 1 & 1   \\  
\end{tabular}
\label{f10b}
}
\caption{a) Truth table for the proposition $((q_1 \to q_2) \land q_1) \to q_2$, and b) Membership table for the set $((C_1 \naturalto C_2) \cap C_1 ) \naturalto C_2$}
\label{f10}
\end{figure}

The law of propositional calculus $((q_1 \to q_2) \land \overline{q}_2) \to \overline{q}_1$ is known as ``modus tollendo tollens''. It can be considered a way of reasoning such that if there is conditional relation $q_1 \to q_2$ and it is negated (\textit{tollendo}) that the consequent $q_2$ is true, it is possible to deduce that the antecedent $q_1$ is false. So, this law negates (\textit{tollens}) that antecedent.

The truth table corresponding to the law of propositional calculus $((q_1 \to q_2) \land \overline{q}_2) \to \overline{q}_1$ is presented in figure \ref{f11a} and the membership table for the set $((C_1 \naturalto C_2) \cap \xvec{C}_2) \naturalto \xvec{C}_1$ which is isomorphic to it and equal to the universal set within whose frame $C_1$ and $C_2$ have been characterized, is presented in figure \ref{f11b}. 

\begin{figure}[H]
\centering
\subfloat[Truth table for the proposition $((q_1 \to q_2) \land \overline{q}_2) \to \overline{q}_1 $]{
\begin{tabular}{c|c|c|c|c|c||c} 
$q_1$ &  $q_2$ & $\overline{q}_1$  & $\overline{q}_2$  & $q_1 \to q_2$  & $ (q_1 \to q_2) \land \overline{q}_2 $ & $((q_1 \to q_2) \land \overline{q}_2) \to \overline{q}_1$\\
\midrule
0 & 0 & 1 & 1 & 1 & 1 & 1 \\  
0 & 1 & 1 & 0 & 1 & 0 & 1   \\  
1 & 0 & 0 & 1 & 0 & 0 & 1   \\  
1 & 1 & 0 & 0 & 1 & 0 & 1   \\  
\end{tabular}
\label{f11a}
}
\\
\subfloat[Membership table for the set $((C_1 \naturalto C_2) \cap \xvec{C}_2) \naturalto \xvec{C}_1$; $(((C_1 \naturalto C_2) \cap \xvec{C}_2) \naturalto \xvec{C}) = \mathbb U$]{
\hspace{-.3in}
\begin{tabular}{c|c|c|c|c|c||c}
$C_1$ &  $C_2$ & $\xvec{C}_1$ & $\xvec{C}_2$ & $C_1 \naturalto C_2$ & $ (C_1 \naturalto C_2) \cap \xvec{C}_2$  & $((C_1 \naturalto C_2) \cap \xvec{C}_2 ) \naturalto \xvec{C}_1$\\
\midrule
0 & 0 & 1 & 1 & 1 & 1 & 1 \\  
0 & 1 & 1 & 0 & 1 & 0 & 1   \\  
1 & 0 & 0 & 1 & 0 & 0 & 1   \\  
1 & 1 & 0 & 0 & 1 & 0 & 1   \\  
\end{tabular}
\label{f11b}%
}
\caption{a) Truth table for the proposition $((q_1 \to q_2) \land \overline{q}_2)\to \overline{q}_1$, and b) Membership table for the set $((C_1 \naturalto C_2) \cap \xvec{C}_2) \naturalto \xvec{C}_1$}
\label{f11}
\end{figure}

Two more examples of laws of propositional calculus are those of both De Morgan's laws. The truth table for one of them is presented in figure \ref{f12a} -- $\overline{(q_1 \lor q_2)} \longleftrightarrow (\overline{q}_1 \land \overline{q}_2)$
-- and the membership table for the set $\xvec{(C_1 \cup C_2)} \naturaltolr (\xvec{C}_1 \cap \xvec{C}_2)$, which is isomorphic to that law and and equal to the universal set within whose frame $C_1$ and $C_2$ have been characterized, is presented in in figure \ref{f12b}.

\begin{figure}[H]
\centering
\subfloat[Truth table for the proposition $\overline{(q_1 \lor q_2)} \longleftrightarrow (\overline{q}_1 \land \overline{q}_2)$]{
\begin{tabular}{c|c|c|c|c|c|c||c|} 
$q_1$ & $q_2$ & $\overline{q}_1$ & $\overline{q}_2$ & $q_1 \lor q_2$ & $\overline{(q_1 \lor q_2)}$ & $\overline{q}_1 \land \overline{q}_2$ & $\overline{(q_1 \lor q_2)} \longleftrightarrow (\overline{q}_1 \land \overline{q}_2)$\\
\midrule
0 & 0 & 1 & 1 & 0 & 1 & 1 & 1  \\  
0 & 1 & 1 & 0 & 1 & 0 & 0 & 1  \\  
1 & 0 & 0 & 1 & 1 & 0 & 0 & 1  \\  
1 & 1 & 0 & 0 & 1 & 0 & 0 & 1   \\  
\end{tabular}
\label{f12a}
}
\\
\subfloat[Membership table for the set $\xvec{(C_1 \cup C_2)} \naturaltolr (\xvec{C}_1 \cap \xvec{C}_2)$; $(\xvec{(C_1 \cup C_2)} \naturaltolr (\xvec{C}_1 \cap \xvec{C}_2))=\mathbb U$ ]{
\hspace{-.1in}
\begin{tabular}{c|c|c|c|c|c|c||c|}
$C_1$ &  $C_2$ & $\xvec{C}_1$ &  $\xvec{C}_2$ & $C_1 \cup C_2$ & $\xvec{(C_1 \cup C_2)}$ & $\xvec{C}_1 \cap \xvec{C}_2$ & $\xvec{(C_1 \cup C_2)} \naturaltolr (\xvec{C}_1 \cap \xvec{C}_2)$\\
\midrule
0 & 0 & 1 & 1 & 0 & 1 & 1 & 1  \\  
0 & 1 & 1 & 0 & 1 & 0 & 0 & 1  \\  
1 & 0 & 0 & 1 & 1 & 0 & 0 & 1  \\  
1 & 1 & 0 & 0 & 1 & 0 & 0 & 1   \\
\end{tabular}
\label{f12b}%
}
\caption{a) Truth table for the proposition $\overline{(q_1 \lor q_2)} \longleftrightarrow (\overline{q}_1 \land \overline{q}_2)$, and b) Membership table for the set $\xvec{(C_1 \cup C_2)} \naturaltolr (\xvec{C}_1 \cap \xvec{C}_2)$}
\label{f12}
\end{figure}

The truth table for the other De Morgan's law -- $\overline{(q_1 \land q_2)} \longleftrightarrow (\overline{q}_1 \land \overline{q}_2)$ -- is presented in figure \ref{f13a}, and the membership table for the set $\xvec{(C_1 \cap C_2)} \naturaltolr (\xvec{C}_1 \cup \xvec{C}_2)$, which is isomorphic to that law and and equal to the universal set within whose frame $C_1$ and $C_2$ have been characterized, is presented in figure \ref{f13b}.

\begin{figure}[H]
\centering
\subfloat[Truth table for the proposition $\overline{(q_1 \land q_2)} \longleftrightarrow (\overline{q}_1 \lor \overline{q}_2)$]{
\begin{tabular}{c|c|c|c|c|c|c||c|} 
$q_1$ & $q_2$ & $\overline{q}_1$ & $\overline{q}_2$ & $q_1 \land q_2$ & $\overline{(q_1 \land q_2)}$ & $\overline{q}_1 \lor \overline{q}_2$ & $\overline{(q_1 \land q_2)} \longleftrightarrow (\overline{q}_1 \lor \overline{q}_2)$\\
\midrule
0 & 0 & 1 & 1 & 0 & 1 & 1 & 1  \\  
0 & 1 & 1 & 0 & 0 & 1 & 1 & 1  \\  
1 & 0 & 0 & 1 & 0 & 1 & 1 & 1  \\  
1 & 1 & 0 & 0 & 1 & 0 & 0 & 1   \\  
\end{tabular}
\label{f13a}
}\\

\subfloat[Membership table for the set $\xvec{(C_1 \cap C_2)} \naturaltolr (\xvec{C}_1 \cup \xvec{C}_2)$; $((\xvec{(C_1 \cap C_2)} \naturaltolr (\xvec{C}_1 \cup \xvec{C}_2)) = \mathbb U$]{
\hspace{-0.3in}
\begin{tabular}{c|c|c|c|c|c|c||c|}
$C_1$ &  $C_2$ & $\xvec{C}_1$ &  $\xvec{C}_2$ & $C_1 \cap C_2$ & $\xvec{(C_1 \cap C_2)}$ & $\xvec{C}_1 \cup \xvec{C}_2$ & $\xvec{(C_1 \cap C_2)} \naturaltolr (\xvec{C}_1 \cup \xvec{C}_2)$\\
\midrule
0 & 0 & 1 & 1 & 0 & 1 & 1 & 1  \\  
0 & 1 & 1 & 0 & 0 & 1 & 1 & 1  \\  
1 & 0 & 0 & 1 & 0 & 1 & 1 & 1  \\  
1 & 1 & 0 & 0 & 1 & 0 & 0 & 1   \\
\end{tabular}

\label{f13b}%
}
\caption{a) Truth table for the proposition $\overline{(q_1 \land q_2)} \longleftrightarrow (\overline{q}_1 \lor \overline{q}_2)$, and b) Membership table for the set $\xvec{(C_1 \cap C_2)} \naturaltolr (\xvec{C}_1 \cup \xvec{C}_2)$}
\label{f13}
\end{figure}

Another law of propositional calculus is the proposition $\overline{(q_1 \dot{\lor} q_2)} \longleftrightarrow (q_1 \longleftrightarrow q_2)$, which expresses the logical equivalence of a) the negation of the exclusive disjunction of the propositions $q_1$ and $q_2$, and b) the logical equivalence of $q_1$ and $q_2$. The truth table for that law is presented in figure \ref{f14a}, and the membership table for the set $\xvec{(C_1 \dot{\cup} C_2)} \naturaltolr (C_1 \naturaltolr C_2)$, which is isomorphic to it and equal to the universal set within whose frame $C_1$ and $C_2$ have been characterized is presented in figure \ref{f14b}.

\begin{figure}[H]
\centering
\subfloat[Truth table for the proposition $\overline{(q_1 \dot{\lor} q_2)} \longleftrightarrow (q_1 \longleftrightarrow q_2)$]{
\begin{tabular}{c|c|c|c|c||c|} 
$q_1$ & $q_2$ & $q_1 \dot{\lor} q_2$ & $\overline{(q_1 \dot{\lor} q_2)}$ & $(q_1 \longleftrightarrow q_2)$ &  $\overline{(q_1 \dot{\lor} q_2)} \longleftrightarrow (q_1 \longleftrightarrow q_2)$\\
\midrule
0 & 0 & 0 & 1 & 1 & 1  \\  
0 & 1 & 1 & 0 & 0 & 1  \\  
1 & 0 & 1 & 0 & 0 & 1  \\  
1 & 1 & 0 & 1 & 1 & 1  \\  
\end{tabular}
\label{f14a}
}
\\
\subfloat[Membership table for the set $\xvec{(C_1 \dot{\cup} C_2)} \naturaltolr (C_1 \naturaltolr C_2)$]{
\hspace{-.1in}
\begin{tabular}{c|c|c|c|c||c|}
$C_1$ &  $C_2$ & $C_1 \dot{\cup} C_2$ & $\xvec{(C_1 \dot{\cup} C_2)}$ & $C_1 \naturaltolr C_2)$ & $\xvec{(C_1 \dot{\cup} C_2)} \naturaltolr (C_1 \naturaltolr C_2)$\\
\midrule
0 & 0 & 0 & 1 & 1 & 1  \\  
0 & 1 & 1 & 0 & 0 & 1  \\  
1 & 0 & 1 & 0 & 0 & 1  \\  
1 & 1 & 0 & 1 & 1 & 1  \\ 
\end{tabular}
\label{f14b}
}
\caption{a) Truth table for the proposition $\overline{(q_1 \dot{\lor} q_2)} \longleftrightarrow (q_1 \longleftrightarrow q_2)$, and b) Membership table for the set $\xvec{(C_1 \dot{\cup} C_2)} \naturaltolr (C_1 \naturaltolr C_2)$}
\label{f14}
\end{figure}

The law of transitivity of material implication in propositional calculus is as follows: $((q_1 \to q_2) \land (q_2 \to q_3)) \to (q_1 \to q_3)$.

The truth table for that proposition is presented in figure \ref{f15a} and the membership table of the corresponding set is presented in figure \ref{f15b}.

\begin{figure}[H]
\centering
\begin{tabular}{c|c|c|c|c|c|}
$q_1$ & $q_2$ & $q_3$ & $q_1 \to q_2$ & $q_2 \to q_3$ & $q_1 \to q_3$\\
\midrule
0 & 0 & 0 & 1 & 1 & 1  \\  
0 & 0 & 1 & 1 & 1 & 1  \\  
0 & 1 & 0 & 1 & 0 & 1   \\  
0 & 1 & 1 & 1 & 1 & 1   \\ 
1 & 0 & 0 & 0 & 1 & 0   \\  
1 & 0 & 1 & 0 & 1 & 1 \\  
1 & 1 & 0 & 1 & 0 & 0 \\  
1 & 1 & 1 & 1 & 1 & 1  \\   
\end{tabular}\\
\subfloat[Truth table for the proposition $((q_1 \to q_2) \land (q_2 \to q_3)) \to (q_1 \to q_3) $; $f(q_1, q_2, q_3) = ((q_1 \to q_2) \land (q_2 \to q_3)) \to (q_1 \to q_3)$]{
\hspace{1in}

\begin{tabular}{c||c|}
 $ (q_1 \to q_2) \land (q_2 \to q_3)$ & $f(q_1, q_2, q_3)$\\
\midrule
  1 & 1\\  
 1 & 1 \\  
 0 & 1 \\  
 1 & 1 \\ 
  0 & 1 \\  
  0 & 1 \\  
  0 & 1 \\  
1 & 1 \\   
\end{tabular}
\label{f15a}
}
\vspace{0.2in}
\\
\begin{tabular}{c|c|c|c|c|}
$C_1$ & $C_2$ & $C_3$ & $C_1 \naturalto C_2$ & $C_2 \naturalto C_3$ \\
\midrule
0 & 0 & 0 & 1 & 1  \\  
0 & 0 & 1 & 1 & 1  \\  
0 & 1 & 0 & 1 & 0 \\  
0 & 1 & 1 & 1 & 1  \\ 
1 & 0 & 0 & 0 & 1 \\  
1 & 0 & 1 & 0 & 1  \\  
1 & 1 & 0 & 1 & 0  \\  
1 & 1 & 1 & 1 & 1 \\   
\end{tabular}\\
\subfloat[Membership table for the set $((C_1 \naturalto C_2) \cap C_2 \naturalto C_3)) \naturalto (C_1 \naturalto C_3)$;  $f(C_1, C_2, C_3)\\= (((C_1 \naturalto C_2) \cap C_2 \naturalto C_3)) \naturalto (C_1 \naturalto C_3)) = \mathbb U$]{
\hspace{1in}
\begin{tabular}{c|c||c|}
 $C_1 \naturalto C_3$ & $ (C_1 \naturalto C_2) \cap (C_2 \naturalto C_3)$ & $f(C_1, C_2, C_3)$\\
\midrule
  1 & 1 & 1\\  
  1 & 1 & 1 \\  
 1 & 0 & 1 \\  
  1 & 1 & 1 \\ 
 0 & 0 & 1 \\  
 1 & 0 & 1 \\  
  0 & 0 & 1 \\  
 1 & 1 & 1 \\   
\end{tabular}
\label{f15b}
}
\caption{a) Truth table for the proposition $((q_1 \to q_2) \land (q_2 \to q_3)) \to (q_1 \to q_3) $, and b) Membership table for the set $((C_1 \naturalto C_2) \cap C_2 \naturalto C_3)) \naturalto (C_1 \naturalto C_3)$}
\label{f15}
\end{figure}

The negation of any law of propositional calculus -- denominated ``contradiction'' -- is a false proposition, regardless of whether the truth values of the propositions of which it is a function. The negation of any contradiction in propositional calculus is a law of that calculus. In effect, the negation of any contradiction in propositional calculus is the negation of the negation (that is, a ``double negation'') of some law of propositional calculus, and is, therefore, equal to that law.

Given that a set corresponding to any law of propositional calculus, isomorphic to it, is equal to the universal set, its complement is the empty set. Just as a set equal to the universal set corresponds to each law of propositional calculus, the empty set corresponds to each contradiction in propositional calculus.

\section{Propositional Calculus of Canonical Fuzzy Logic (CFL)}

There is literature such as that cited above in which reasons are given for why in many cases it is useful to assign, to each proposition, not only one of two possible truth values (true or false) but rather a degree (or weight) of truth that, as accepted, varies between 0 and 1. This is one of the basic notions of fuzzy logic: The weight of truth of the proposition $q$ will be symbolized as $w(q)$.

Consider the ``extreme'' cases of $w(q)$: $w(q)=0$ and $w(q)=1$. With $w(q)=0$, it is known that $q$ is considered absolutely not true; that is, totally false. With $w(q)=1$ it is known that $q$ is considered absolutely true.

It will be accepted that $w(\overline{q}) + w(q) = 1$. Therefore, $w(q) =1-w(\overline{q})$ and $w(\overline{q})=1-w(q)$.

In the table in figure \ref{f16} and in the text that follows it, an explanation is provided about how to assign a weight of truth to the proposition which is the inclusive disjunction of the propositions $q_1$ and $q_2$ (that is, $q_1 \lor q_2$, according to CFL). That weight of truth is symbolized as $w(q_1 \lor q_2)$.
\begin{figure}[H]
\centering
\hspace{0.8in}
\begin{tabular}{c|c||c|l}
$q_1$ & $q_2$ & $q_1\lor q_2$ & $w(q_1 \lor q_2) = \displaystyle\sum_{i=1}^{4}S_i$\\
\midrule
0 & 0 & 0 & $S_1 = 0$  \\
0 & 1 & 1 & $S_2 = w(\overline{q}_1) \cdot w(q_2)$  \\  
1 & 0 & 1 & $S_3 = w(q_1) \cdot w(\overline{q}_2)$  \\  
1 & 1  & 1 & $S_4 = w(q_1) \cdot w(q_2)$\\  
\end{tabular}
\hspace{1in}
\caption{Truth table for $q_1 \lor q_2$ and the corresponding $w(q_1 \lor q_2)$, according to CFL}
\label{f16}

\end{figure}

Note that the first three columns of the table in figure \ref{f16} make up the truth table for $q_1 \lor q_2$ according to classical bivalent logic (BL). In addition, in each of the other tables considered in this section, the first three columns make up the truth table according to classical bivalent logic (BL) of the proposition to which reference is made. 

The addends -- or summands -- equal to 0 correspond to the numerical values of 0 in the column with the proposition whose weight of truth according to CFL is to be computed -- in this case, $q_1 \lor q_2$. Here the only addend of this type is $S_1$:
$S_1=0$.

Addends, which here are $S_2$,  $S_3$, and $S_4$ -- computed as explained below -- correspond to the numerical values of 1 in the column with the proposition whose weight of truth according to CFL is to be computed (in this case, $q_1 \lor q_2$).

Recall that a numerical value of 0 in the column corresponding to $q_j$, for $j=1, 2$, specifies the negation of $q_j$ -- that is, $\overline{q}_j$ -- and the weight of truth of this proposition is symbolized as $w(\overline{q}_j)$. Recall also that a numerical value of 1 in the column corresponding to $q_j$, for $j=1, 2$, specifies the affirmation of 
$q_j$ and the weight of truth of this proposition is symbolized as $w(q_j)$. Therefore, if a 0 is assigned to $q_1$, and a 1 is assigned to $q_2$, as in the second row of the table in figure 16, the corresponding weights of truth are $w(\overline{q}_1)$ and $w(q_2)$. According to CFL the numerical value of $S_2$ is the product of those weights of truth: $S_2=w(\overline{q}_1) \cdot w(q_2)$.

If a 1 is assigned to $q_1$, and a 0 is assigned to $q_2$, as in the third row of the table in figure 16, the corresponding weights of truth are $w(q_1)$ and $w(\overline{q}_2)$. According to CFL the numerical value of $S_3$ is the product of those weights of truth: $S_3=w(q_1) \cdot w(\overline{q}_2)$.

If a 1 is assigned to $q_1$, and a 1 is also assigned to $q_2$, as occurs here in the fourth row of the table in figure 16, the corresponding weights of truth are $w(q_1)$ and $w(_2)$. According to CFL the numerical value of $S_4$ is the product of those weights of truth: $S_4=w(q_1) \cdot w(q_2)$.

Given $S_1$, $S_2$, $S_3$ and $S_4$, $w(q_1 \lor w_2)$ can be computed:\\
\noindent $w(q_1 \lor w_2) = \displaystyle\sum_{i=1}^{4}S_i = S_1 + S_2 + S_3 + S_4 = 0 + w(\overline{q}_1) \cdot w(q_2) + w(q_1) \cdot w(\overline{q}_2) + w(q_1) \cdot w(q_2) = (1-w(q_1)) \cdot w(q_2) + w(q_1) \cdot (1-w(q_2)) + w(q_1) \cdot w(q_2) = w(q_2)-w(q_1)\cdot w(q_2) + w(q_1)-w(q_1)\cdot w(q_2) + w(q_1) \cdot w(q_2) = w(q_1) + w(q_2) - 2 \cdot (w(q_1)\cdot w(q_2)) + w(q_1)\cdot w(q_2) = w(q_1) + w(q_2) - w(q_1)\cdot w(q_2)$\\

\noindent That is, according to CFL, the value of truth of $w(q_1 \lor q_2)$ is the following:  $w(q_1\lor q_2) = w(q_1) + w(q_2) - w(q_1)\cdot w(q_2)$.

Note that the above equation is different from that accepted as the currently most widespread equation in fuzzy logic \cite{d5}: $w(q_1 \lor q_2) =$ Max $\{ w(q_1), w(q_2) \}$.

Suppose, for example, that it is accepted that for certain cases of $q_1$ and $q_2$, $w(q_1) = 0.8$ and $w(q_2)=0.6$. Then, $w(q_1\lor q_2)=0.8 + 0.6 - (0.8)\cdot(0.6) = 0.92$, according to CFL, whereas $w(q_1 \lor q_2) = $ Max $\{ 0.8, 0.6 \} = 0.8$, according to the currently best known version of fuzzy logic.

In the table in figure \ref{f17} it is shown how to assign a weight of truth, according to CFL, to the proposition which is the conjunction of the two propositions $q_1$ and $q_2$; that is, $q_1 \land q_2$.

\begin{figure}[H]
\centering
\hspace{1in}
\begin{tabular}{c|c||c|l}
$q_1$ & $q_2$ & $q_1\land q_2$ & $w(q_1 \land q_2) = \displaystyle\sum_{i=1}^{4}S_i$\\
\midrule
0 & 0 & 0 & $S_1=0$  \\
0 & 1 & 1 & $S_2=0$  \\  
1 & 0 & 1 & $S_3=0$  \\  
1 & 1  & 1 & $S_4=w(q_1) \cdot w(q_2)$\\  
\end{tabular}
\hspace{1in}
\caption{Truth table for $q_1 \land q_2$ and the corresponding $w(q_1 \land q_2)$, according to CFL}
\label{f17}

\end{figure}
Note that using the same criterion as that of the table above (that of assigning the numerical value of 0 to each addend $S_i$, for $i=1, 2, 3, 4$, such that there is a 0 in the column corresponding to the proposition whose weight of truth is to be computed, and in the same row of the addend considered), the summands $S_1$, $S_2$ and $S_3$ are equal to 0: $S_1= S_2=S_3=0$. Given that in that column (that of $q_1 \land q_2$) and in that same row corresponding to the addend $S_4$ there is a 1, $S_4$ is equal to the product of the weights of truth of $q_1$ and $q_2$: $S_4 = w(q_1) \cdot w(q_2)$. Hence, $w(q_1 \land q_2)$ can be computed as follows:

\noindent $w(q_1 \land q_2) = \displaystyle\sum_{i=1}^{4}S_i = 0 + 0 + 0 + (w(q_1) \cdot w(q_2))$\\

\noindent Therefore, $w(q_1 \land q_2) = w(q_1) \cdot w(q_2)$, according to CFL.

Note that the above equation is different from the currently best known version in fuzzy logic \cite{d5}: $w(q_1 \land q_2) =$ Min $ \{w(q_1), w(q_2) \}$.
Suppose, for example, that it is accepted that for certain cases of $q_1$ and $q_2$, $w(q_1)=0.8$ and $w(q_2)=0.6$. Then, $w(q_1 \land q_2) = w(q_1) \cdot w(q_2) = (0.8) \cdot (0.6) = 0.48$, according to CFL; whereas $w(q_1 \land q_2) = $ Min $ \{w(q_1), w(q_2) \}$ = Min $ \{(0.8), (0.6)\} = 0.6$, in the best known version of fuzzy logic.

Consideration will be given below in this section to some propositions resulting from other operations of propositional calculus, and their corresponding truth values will be computed, according to CFL, with the same method used for $q_1 \lor q_2$ and $q_1 \land q_2$.

In the table in figure \ref{f18}, it is shown how to assign a weight of truth, according to CFL, to the proposition which is the exclusive disjunction of the propositions  $q_1$ and $q_2$ (that is, $q_1 \dot{\lor} q_2$). 
\begin{figure}[H]
\centering
\hspace{1in}
\begin{tabular}{c|c||c|l}
$q_1$ & $q_2$ & $q_1 \dot{\lor} q_2$ & $w(q_1 \dot{\lor} q_2) = \displaystyle\sum_{i=1}^{4}S_i$\\
\midrule
0 & 0 & 0 & $S_1=0$  \\
0 & 1 & 1 & $S_2=w(\overline{q}_1) \cdot w(q_2)$\\   
1 & 0 & 1 & $S_3=w(q_1) \cdot w(\overline{q}_2)$\\  
1 & 1 & 0 & $S_4=0$ 
\end{tabular}
\hspace{1in}
\caption{Truth table for $q_1 \dot{\lor} q_2$ and the corresponding $w(q_1 \dot{\lor} q_2)$, according to CFL}
\label{f18}

\end{figure}

Given $S_1, S_2, S_3$ and $S_4$, $w(q_1 \dot{\lor} q_2)$ can be computed according to CFL:

\noindent $w(q_1 \dot{\lor} q_2)=\displaystyle\sum_{i=1}^{4}S_i=S_1 + S_2 + S_3 + S_4 = 0 + w(\overline{q}_1) \cdot w(q_2) + w(q_1) \cdot w(\overline{q}_2) + 0 = (1-w(q_1)) \cdot w(q_2) + w(q_1) \cdot (1-w(q_2)) = w(q_2)-w(q_1)\cdot w(q_2) + w(q_1)-w(q_1)\cdot w(q_2) = w(q_1) + w(q_2) - 2 \cdot w(q_1)\cdot w(q_2)$\\

\noindent That is, $w(q_1 \dot{\lor} q_2) = w(q_1) + w(q_2) - 2 \cdot w(q_1)\cdot w(q_2)$.

In the table in figure \ref{f19} it is shown how to assign a weight of truth to the conditional proposition $q_1 \to q_2$, in which $q_1$ is the antecedent of that material implication and $q_2$ is its consequent.

\begin{figure}[H]
\centering
\hspace{1in}
\begin{tabular}{c|c||c|l}
$q_1$ & $q_2$ & $q_1 \to q_2$ & $w(q_1 \to q_2) = \displaystyle\sum_{i=1}^{4}S_i$\\
\midrule
0 & 0 & 1 & $S_1=w(\overline{q}_1) \cdot w(\overline{q}_2)$  \\
0 & 1 & 1 & $S_2=w(\overline{q}_1) \cdot w(q_2)$\\   
1 & 0 & 0 & $S_3=0$\\  
1 & 1 & 1 & $S_4=w(q_1) \cdot w(q_2)$ 
\end{tabular}
\hspace{1in}
\caption{Truth table for $q_1 \to q_2$ and the corresponding $w(q_1 \to q_2)$, according to CFL}
\label{f19}

\end{figure}
Given $S_1, S_2, S_3$ and $S_4$, $w(q_1) \to q_2)$ can be computed according to CFL:
\noindent $w(q_1 \to q_2)=\displaystyle\sum_{i=1}^{4}S_i=S_1 + S_2 + S_3 + S_4 = w(\overline{q}_1) \cdot w(\overline{q}_2) + w(\overline{q}_1) \cdot w(q_2) + 0 + w(q_1) \cdot w(q_2) = (1-w(q_1)) \cdot (1-w(q_2)) + (1 - w(q_1)) \cdot w(q_2) + w(q_1) \cdot w(q_2) = 1-w(q_2)-w(q_1)+ w(q_1)\cdot w(q_2) + w(q_2)-w(q_1)\cdot w(q_2) + w(q_1) \cdot (w(q_2) = 1-w(q_1) + w(q_1) \cdot (w(q_2)$\\

\noindent That is, $w(q_1 \to q_2) = 1 - w(q_1) + w(q_1) \cdot w(q_2)$.

In the table in figure \ref{f20} it is shown how to assign a weight of truth to the conditional proposition $q_2 \to q_1$, in which $q_2$ is the antecedent of that material implication and $q_1$ is its consequent.

\begin{figure}[H]
\centering
\hspace{0.8in}
\begin{tabular}{c|c||c|l}
$q_1$ & $q_2$ & $q_2 \to q_1$ & $w(q_2 \to q_1) = \displaystyle\sum_{i=1}^{4}S_i$\\
\midrule
0 & 0 & 1 & $S_1=w(\overline{q}_1) \cdot w(\overline{q}_2)$  \\
0 & 1 & 0 & $S_2=0$\\   
1 & 0 & 1 & $S_3=w(q_1) \cdot w(\overline{q}_2)$\\  
1 & 1 & 1 & $S_4=w(q_1) \cdot w(q_2)$ 
\end{tabular}
\hspace{1in}
\caption{Truth table for $q_2 \to q_1$ and the corresponding $w(q_2 \to q_1)$, according to CFL}
\label{f20}

\end{figure}

Given $S_1, S_2, S_3$ and $S_4$, $w(q_2 \to q_1)$ can be computed according to CFL:

\noindent $w(q_2 \to q_1)=\displaystyle\sum_{i=1}^{4}S_i=S_1 + S_2 + S_3 + S_4 = w(\overline{q}_1) \cdot w(\overline{q}_2) + 0 + w(q_1) \cdot w(\overline{q}_2) + w(q_1) \cdot w(q_2) = (1-w(q_1)) \cdot (1-w(q_2)) + w(q_1) \cdot (1-w(q_2)) + w(q_1) \cdot w(q_2) = 1-w(q_2)-w(q_1)+ w(q_1)\cdot w(q_2) + w(q_1)-w(q_1)\cdot w(q_2) + w(q_1) \cdot (w(q_2) = 1-w(q_2) + w(q_1) \cdot (w(q_2)$\\

\noindent That is, $w(q_2 \to q_1) = 1 - w(q_2) + w(q_1) \cdot w(q_2)$.

In the table in figure \ref{f21} it is shown how to assign a weight of truth to the bi-conditional proposition $q_1 \longleftrightarrow q_2$; that is, to the logical equivalence (or material bi-implication) of $q_1$ and $q_2$.
\begin{figure}[H]
\centering
\hspace{0.7in}
\begin{tabular}{c|c||c|l}
$q_1$ & $q_2$ & $q_1 \longleftrightarrow q_2$ & $w(q_1 \longleftrightarrow q_2) = \displaystyle\sum_{i=1}^{4}S_i$\\
\midrule
0 & 0 & 1 & $S_1=w(\overline{q}_1) \cdot w(\overline{q}_2)$  \\
0 & 1 & 0 & $S_2=0$\\   
1 & 0 & 0 & $S_3=0$\\  
1 & 1 & 1 & $S_4=w(q_1) \cdot w(q_2)$ 
\end{tabular}
\hspace{1in}
\caption{Truth table for $q_1 \to q_2$ and the corresponding $w(q_1 \to q_2)$, according to CFL}
\label{f21}

\end{figure}

Given $S_1, S_2, S_3$ and $S_4$, $w(q_1 \longleftrightarrow q_2)$ can be computed according to CFL:

\noindent $w(q_1 \longleftrightarrow q_2)=\displaystyle\sum_{i=1}^{4}S_i=S_1 + S_2 + S_3 + S_4 = w(\overline{q}_1) \cdot w(\overline{q}_2) + 0 + 0 + w(q_1) \cdot w(q_2) = (1-w(q_1)) \cdot (1 - w(q_2)) + w(q_1) \cdot w(q_2) = 1-w(q_2)-w(q_1) + w(q_1)\cdot w(q_2) + w(q_1)\cdot w(q_2) = 1-w(q_1) - w(q_2) + 2\cdot (w(q_1)\cdot w(q_2)$\\

\noindent That is, \noindent $w(q_1 \longleftrightarrow q_2) = 1 - w(q_1) - w(q_2) + 2\cdot w(q_1)\cdot w(q_2)$.

In the table in figure \ref{f22} it is shown how to assign a weight of truth to the proposition $q_1$ nand $q_2$ -- that is,  $q_1 \uparrow q_2$. Recall that the operator symbolized as $ \uparrow$ is also known as \textit{Sheffer stroke}.

\begin{figure}[H]
\centering
\hspace{1in}
\begin{tabular}{c|c||c|l}
$q_1$ & $q_2$ & $q_1 \uparrow q_2$ & $w(q_1 \uparrow q_2) = \displaystyle\sum_{i=1}^{4}S_i$\\
\midrule
0 & 0 & 1 & $S_1=w(\overline{q}_1) \cdot w(\overline{q}_2)$\\
0 & 1 & 1 & $S_2=w(\overline{q}_1) \cdot (q_2)$\\   
1 & 0 & 1 & $S_3=w(q_1) \cdot w(\overline{q}_2)$\\  
1 & 1 & 0 & $S_4=0$ 
\end{tabular}
\hspace{1in}
\caption{Truth table for $q_1 \uparrow q_2$ and the corresponding $w(q_1 \uparrow q_2)$, according to CFL}
\label{f22}

\end{figure}

Given $S_1,  S_2, S_3$ and $S_4$, $w(q_1 \uparrow q_2)$ can be computed according to CFL:

\noindent  $w(q_1 \uparrow q_2)=\displaystyle\sum_{i=1}^{4}S_i=S_1 + S_2 + S_3 + S_4 = w(\overline{q}_1) \cdot w(\overline{q}_2) + w(\overline{q}_1)\cdot w(q_2) + w(q_1) \cdot w(\overline{q}_2) + 0 = (1-w(q_1)) \cdot (1 - w(q_2)) + (1-w(q_1)) \cdot w(q_2) + w(q_1) \cdot (1-w(q_2)) = 1-w(q_2) - w(q_1) + w(q_1)\cdot w(q_2) + w(q_2) - w(q_1) \cdot w(q_2) + w(q_1) - w(q_1) \cdot w(q_2) = 1- w(q_1)\cdot w(q)_2$\\

\noindent That is, $w(q_1 \uparrow q_2) = 1 - w(q_1)\cdot w(q_2)$.

In the table in figure \ref{f23} it is shown how to assign a weight of truth to the proposition $q_1$ nor $q_2$ -- that is, $q_1 \downarrow q_2$. Recall that the operator symbolized as $ \downarrow$ is also known as \textit{Peirce's arrow}.

\begin{figure}[H]
\centering
\hspace{1in}
\begin{tabular}{c|c||c|l}
$q_1$ & $q_2$ & $q_1 \downarrow q_2$ & $w(q_1 \downarrow q_2) = \displaystyle\sum_{i=1}^{4}S_i$\\
\midrule
0 & 0 & 1 & $S_1=w(\overline{q}_1) \cdot w(\overline{q}_2)$\\
0 & 1 & 0 & $S_2=0$\\   
1 & 0 & 0 & $S_3=0$\\  
1 & 1 & 0 & $S_4=0$ 
\end{tabular}
\hspace{1in}
\caption{Truth table for $q_1 \downarrow q_2$ and the corresponding weight of truth $w(q_1 \downarrow q_2)$, according to CFL}
\label{f23}

\end{figure}

Given $S_1, S_2, S_3$ and $S_4$, $w(q_1 \downarrow q_2)$ can be computed according to CFL:

\noindent $w(q_1 \downarrow q_2)=\displaystyle\sum_{i=1}^{4}S_i=S_1 + S_2 + S_3 + S_4 = w(\overline{q}_1) \cdot w(\overline{q}_2) + 0 + 0 +0 = (1-w(q_1)) \cdot (1 - w(q_2)) =1- w(q_2) - w(q_1) + w(q_1)\cdot w(q_2)$\\

\noindent That is, $w(q_1 \downarrow q_2) = 1 - w(q_1) - w(q_2) + w(q_1) \cdot w(q_2)$.

\section{The Laws -- or Tautologies -- of Propositional Calculus of BL Are Also Valid in the Propositional Calculus of CFL}

The laws or tautologies of propositional calculus of BL are also valid in the propositional calculus of CFL according to the following criterion: the weight of truth of each law in the latter calculus is equal to 1.

Consideration is given below to several laws of propositional calculus of CFL.

In the table in figure \ref{f24} it is shown how to assign a weight of truth, according to CFL, to the law of the excluded middle, $q_1 \lor \overline{q}_2$. This law is a function of a sole proposition $q$.
\begin{figure}[H]
\centering
\hspace{0.8in}
\begin{tabular}{c|c||c|l}
$q$ & $\overline{q}$ & $q\lor \overline{q}$ & $w(q \lor \overline{q}) = \displaystyle\sum_{i=1}^{2}S_i$\\
\midrule
0 & 1& 1 & $S_1=w(\overline{q})$  \\  
1 & 0  & 1 & $S_2=w(q) $\\  
\end{tabular}
\hspace{1in}
\caption{Truth table for $q\lor \overline{q}$ and the corresponding weight of truth according to CFL}
\label{f24}
\end{figure}

Since $q\lor \overline{q}$ is a proposition that is a function of a sole proposition $q$, the truth table of this law has only 2 rows: $2^1=2$,

Given $S_1$ and $S_2$, $w(q \lor \overline{q})$ can be computed, according to CFL: $w(q \lor \overline{q}) = \displaystyle\sum_{i=1}^{2}S_i =  S_1 + S_2 = w(\overline{q}) + w(q) = (1 - w(q)) + w(q) = 1 - w(q) + w(q) = 1$.\\

\noindent That is, $w(q \lor \overline{q}) = 1$.

In the table in figure \ref{f25} it is shown how to assign a weight of truth, according to CFL, to the law of non-contradiction $(\overline{q \land \overline{q}})$. This law also is a proposition that is a function of a sole proposition $q$.

\begin{figure}[H]
\centering
\hspace{0.5in}
\begin{tabular}{c|c|c||c|l}
$q$ & $\overline{q}$ & $q\land \overline{q}$ & $(\overline{q \land \overline{q}})$ & $w(\overline{q \land \overline {q}}) = \displaystyle\sum_{i=1}^{2}S_i$\\
\midrule
0 & 1 & 0 & 1 & $S_1=w(\overline{q})$  \\  
1 & 0 & 0 & 1 & $S_2=w(q)$\\  
\end{tabular}
\hspace{1in}
\caption{Truth table for $(\overline{q \land \overline{q}})$ and the corresponding weight of truth}
\label{f25}
\end{figure}

Given $S_1$ and $S_2$, $w(\overline{q \land \overline {q}})$ can be computed, according to CFL:

\noindent $w(\overline{q \land \overline {q}}) = \displaystyle\sum_{i=1}^{2}S_i = S_1 + S_2= w(\overline{q}) + w(q) = 1$\\

\noindent That is, $w(\overline{q \land \overline{q}}) = 1$.

In the table in figure \ref{f26} it is shown how to assign a weight of truth, according to CFL, to one of the De Morgan's laws $(q_1 \lor q_2) \longleftrightarrow (\overline{q}_1 \land \overline{q}_2)$. This law is a proposition which is a function of two propositions, $q_1$ and $q_2$. Therefore, the corresponding truth table has 4 rows $(2^2=4)$, and each addend $S_i$, for $i = 1, 2, 3, 4$, will be a product of two functions.

\begin{figure}[H]
\hspace{0.5in}
\begin{tabular}{c|c|c|c|c|c|c||}
$q_1$ & $q_2$ & $q_1 \lor q_2$ & $(\overline{q_1 \lor q_2})$ & $\overline{q}_1$ & $ \overline {q}_2$ & $(\overline{q}_1 \land \overline {q}_2)$\\
\midrule
0 & 0 & 0 & 1 & 1 & 1 & 1 \\  
0 & 1 & 1 & 0 & 1 & 0 & 0 \\  
1 & 0 & 1 & 0 & 0 & 1 & 0 \\  
1 & 1 & 1 & 0 & 0 & 0 & 0 \\   
\end{tabular}
\end{figure}
\begin{figure}[H]
\hspace{0.2in}
\begin{tabular}{c|l}
$(\overline{q_1 \lor q_2}) \longleftrightarrow (\overline{q}_1 \land \overline {q}_2)$ & $w((\overline{q_1 \lor q_2}) \longleftrightarrow (\overline{q}_1 \land \overline {q}_2)) = \displaystyle\sum_{i=1}^{4}S_i$\\
\midrule
1 & $S_1=w(\overline{q}_1) \cdot w(\overline{q}_2)$\\  
1 & $S_2=w(\overline{q}_1) \cdot w(q_2)$\\  
1 & $S_3=w(q_1) \cdot w(\overline{q}_2)$\\  
1 & $S_4=w(q_1) \cdot w(q_2)$\\   
\end{tabular}
\caption{Truth table for $(\overline{q_1 \lor q_2}) \longleftrightarrow (\overline{q}_1 \land \overline {q}_2)$, and the corresponding weight of truth, according to CFL}
\label{f26}
\end{figure}
Given $S_1, S_2, S_3$ and $S_4$, $w((\overline{q_1 \lor q_2}) \leftrightarrow (\overline{q}_1 \land \overline {q}_2))$ can be computed according to CFL: $w((\overline{q_1 \lor q_2}) \longleftrightarrow (\overline{q}_1 \land \overline{q}_2)) = \displaystyle\sum_{i=1}^{4}S_i = S_1 + S_2 + S_3$ + $S_4 = w(\overline{q}_1) \cdot w(\overline{q}_2) + w(\overline{q}_1) \cdot w(\overline{q}_2) + w(q_1) \cdot w(\overline{q}_2) + w(q_1) \cdot w(q_2) = w(\overline{q}_1) \cdot (w(\overline{q}_2) + w(q_2)) + w(q_1)\cdot(w(\overline{q}_2) + w(q_2)) = w(\overline{q}_1) \cdot 1 + w(q_1) \cdot 1 = w(\overline{q}_1) + w(q_1) = 1$.

In the table in figure \ref{f27} it is shown how to assign the weight of truth, according to CFL, to the law of transitivity of material implication -- $ ((q_1 \to q_2) \land (q_2 \to q_3)) \to (q_1 \to q_3)$. This law is a proposition which is a function of 3 propositions ($q_1, q_2$ and $q_3$). Thus, the corresponding truth table has 8 rows: $2^3 = 8$.
\begin{figure}[H]
\begin{tabular}{c|c|c|c|c|c|c||}
$q_1$ & $q_2$ &   $q_3$ &  $q_1 \to q_2$ &   $q_2 \to q_3$ & $(q_1 \to q_2) \land (q_2 \to q_3)$ & $q_1 \to q_3$\\
\midrule
0 & 0 & 0 & 1 & 1 & 1 & 1 \\  
0 & 0 & 1 & 1 & 1 & 1 & 1  \\ 
0 & 1 & 0 & 1 & 0 & 0 & 1 \\  
0 & 1 & 1 & 1 & 1 & 1 & 1  \\ 
1 & 0 & 0 & 0 & 1 & 0 & 0 \\  
1 & 0 & 1 & 0 & 1 & 0 & 1  \\ 
1 & 1 & 0 & 1 & 0 & 0 & 0 \\  
1 & 1 & 1 & 1 & 1 & 1 & 1  \\ 
\end{tabular}
\end{figure}
\begin{figure}[H]
\hspace{-0.5in}
\begin{tabular}{c|l}
$((q_1 \to q_2) \land (q_2 \to q_3)) \to (q_1 \to q_3)$ & $w((q_1 \to q_2) \land (q_2 \to q_3)) \to (q_1 \to q_3) = \displaystyle\sum_{i=1}^{8}S_i$\\
\midrule
1 & $S_1=w(\overline{q}_1) \cdot w(\overline{q}_2) \cdot w(\overline{q}_3) $\\  
1 & $S_1=w(\overline{q}_1) \cdot w(\overline{q}_2) \cdot w(q_3) $\\  
1 & $S_1=w(\overline{q}_1) \cdot w({q}_2) \cdot w(\overline{q}_3) $\\  
1 & $S_1=w(\overline{q}_1) \cdot w({q}_2) \cdot w({q}_3) $\\  
1 & $S_1=w({q}_1) \cdot w(\overline{q}_2) \cdot w(\overline{q}_3) $\\  
1 & $S_1=w({q}_1) \cdot w(\overline{q}_2) \cdot w({q}_3) $\\  
1 & $S_1=w({q}_1) \cdot w({q}_2) \cdot w(\overline{q}_3) $\\  
1 & $S_1=w({q}_1) \cdot w({q}_2) \cdot w({q}_3) $\\  
\end{tabular}
\caption{Truth table for $((q_1 \to q_2) \land (q_2 \to q_3)) \to (q_1 \to q_3)$, and the corresponding weight of truth, according to CFL}
\label{f27}
\end{figure}

Given $S_1, S_2, S_3, S_4, S_5, S_6, S_7$ and $S_8$, $w(((q_1 \to q_2) \land (q_2 \to q_3)) \to (q_1 \to q_3))$ can be computed according to CFL:

\noindent $w(((q_1 \to q_2) \land (q_2 \to q_3)) \to (q_1 \to q_3)) = \displaystyle\sum_{i=1}^{8}S_i = S_1 + S_2 + S_3 + S_4 + S_5 + S_6 + S_7 + S_8 = w(\overline{q}_1) \cdot w(\overline{q}_2) \cdot w(\overline{q}_3) + w(\overline{q}_1) \cdot w(\overline{q}_2) \cdot w(q_3) + w(\overline{q}_1) \cdot w(q_2) \cdot w(\overline{q}_3) + w(\overline{q}_1) \cdot w(q_2) \cdot w(q_3) + w(q_1) \cdot w(\overline{q}_2) \cdot w(\overline{q}_3) + w(q_1) \cdot w(\overline{q}_2) \cdot w(q_3) + w(q_1) \cdot w(q_2) \cdot w(\overline{q}_3) + w(q_1) \cdot w(q_2) \cdot w(q_3) = w(\overline{q}_1) \cdot (w(\overline{q}_2) \cdot w(\overline{q}_3) + w(\overline{q}_2) \cdot w(q_3) + w(q_2) \cdot w(\overline{q}_3) + w(q_2) \cdot w(q_3))+ w(q_1) \cdot (w(\overline{q}_2) \cdot w(\overline{q}_3) + w(\overline{q}_2) \cdot w(q_3) + w(q_2) \cdot w(\overline{q}_3) + w(q_2) \cdot w(q_3)) = (w(\overline{q}_1) + w(q_1)) \cdot (w(\overline{q}_2) \cdot w(\overline{q}_3) + w(\overline {q}_2) \cdot w(q_3) + w(q_2) \cdot w(\overline{q}_3) + w(q_2) \cdot w(q_3)) = 1 \cdot (w(\overline{q}_2) \cdot (w(\overline{q}_3) + w(q_3)) + w(q_2) \cdot (w(\overline{q}_3) + w(q_3)) = 1 \cdot ((w(\overline{q}_2) \cdot 1) + (w(q_2) \cdot 1))= w(\overline{q}_2) + w(q_2) = 1$\\

In the following considerations, any proposition which is a law of propositional calculus depending on a sole proposition will be symbolized as $l(q_1)$. (The letter $l$ is used for ${l}$aw.) Note that in this case despite considering any law of propositional calculus that is a function of a sole proposition, due to a recursive-type equation which will be presented, the symbol $q_1$ was chosen rather than $q$.

Any law of propositional calculus which is a function of 2 propositions will be symbolized as: $l(q_1, q_2)$. Any law of that calculus which is a function of 3 propositions will be symbolized as: $l(q_1, q_2, q_3)$. Any law of that calculus which is a function of $n$ propositions -- n = 1, 2, 3, \dots -- will be symbolized as: $l(q_1, q_2, q_3, \dots, q_n)$.

According to CFL, the weights of truth of the laws of propositional calculus can be expressed not only as summations of certain addends, but also as products of certain factors. Thus, for example, the weight of truth of any law of propositional calculus which is a function of 2 propositions $q_1$ and $q_2$ can be expressed as: $$ w(l(q_1, q_2)) = \prod_{j=1}^2 (w(\overline{q}_j) + w(q_j))= (w(\overline{q}_1) + w(q_1)) \cdot (w(\overline{q}_2) + w(q_2)) = 1$$

In the multiplication above, each of the factors is equal to 1, so their product is also equal to 1.

The weight of truth of any law of propositional calculus which is a function of 3 propositions $q_1,q_2$ and $q_3$ can be expressed as: $$ w(l(q_1, q_2, q_3)) = \prod_{j=1}^3 (w(\overline{q}_j) + w(q_j)) = (w(\overline{q}_1) + w(q_1)) \cdot (w(\overline{q}_2) + w(q_2)) \cdot (w(\overline{q}_3) + $$ $w(q_3)) = 1$.

In the above multiplication, each factor is equal to 1, so their product is also equal to 1.

In general, according to CFL, the weight of truth of any law of propositional calculus which is the function of $n$ propositions -- $n=1, 2, 3, \dots $ -- can be expressed as follows:

$$ w(l(q_1,q_2, q_3 \dots, q_n)) = \prod_{j=1}^n (w(\overline{q}_j) + w(q_j))$$ $= (w(\overline{q}_1) + w(q_1)) \cdot (w(\overline{q}_2) + w(q_2)) \cdot (w(\overline{q}_3) + w(q_3)) \cdot \dots (w(\overline{q}_n) + w(q_n)) = 1$\\

In the above multiplication, each factor is equal to 1, so their product is also equal to 1.

Note the validity of the following recursive-type equation:

\noindent $w(l(q_1, q_2, q_3, \dots, q_n)) = w(l(q_1, q_2, q_3, \dots, q_{n-1})) \cdot (w(\overline{q_n}) + w(q_n)) = 1$

\noindent That is: 
$$ \prod_{j=1}^n (w(\overline{q}_j) + w(q_j))= \prod_{j=1}^{n-1} (w(\overline{q}_j) + w(q_j)) \cdot (w(\overline{q}_n) + w(q_n)) = 1$$

\section{Fuzzy Sets according to CFL}

\subsection{Propositions Corresponding to Each Element of a Set}

One way to justify, and from a pedagogical standpoint, to facilitate operations with fuzzy sets according to CFL is by considering that each element belonging to a set corresponds to a particular proposition. Recall, in the theory of classical sets, the following way of symbolizing a given set $C_1$: $C_1 = \{x_1, x_2, x_3\}$.

A possible explanation of the meaning of the above equality is as follows: Only the elements $x_1$, $x_2$, and $x_3$ belong to the set $C_1$. Using a familiar notation: $x_1 \in C_1$ ($x_1$ belongs to $C_1$); $x_2 \in C_1$ ($x_2$ belongs to $C_1$); and $x_3 \in C_1$ ($x_3$ belongs to $C_1$). Note that  $x_1 \in C_1$, $x_2 \in C_1$, and $x_3 \in C_1$ are propositions. These three propositions suffice to characterize $C_1$. Given the characterization of $C_1$, it is admitted that each of these three propositions is true.

The proposition $x_1 \in C_1$ can be symbolized in an abbreviated way as $q_{1, 1}$. The propositions $x_2 \in C_1$ and $x_3 \in C_1$ can be symbolized in an abbreviated way as $q_{2, 1}$ and $q_{3,1}$ respectively. In each of these abbreviated symbolizations, the first subscript specifies the element referred to -- 1 refers to $x_1$, 2 to $x_2$ and 3 to $x_3$; and the second subscript specifies the set to which those elements belong: $C_1$.

Thus, for example, the abbreviated symbolizations $q_{2, 17}$, $q_{48, 5}$, and $q_{7, 10}$ mean $x_2 \in C_{17}$, $x_{48} \in C_5$, and $x_7 \in C_{10}$, respectively.

This type of notation will be useful when considering fuzzy sets according to CFL.

\subsection{Fuzzy Finite Sets, according to CFL}

According to CFL, a certain weight of truth (a real number between 0 and 1) should be assigned to each proposition that specifies the membership of a given element in some set.  Thus, for example, a weight of truth such as $0.57$ -- $w(x_2 \in C_1) = 0.57$; abbreviated as $w(q_{2, 1}) = 0. 57$ -- should be assigned to the proposition referring to the element $x_2$, which belongs to the set $C_2$. Recall that $w$ is used for \textit{w}eight.  The weight of truth assigned to each proposition of this type should be determined according to the criterion of experts on the topic considered, except for the propositions related to the elements belonging to two specific sets. This will be considered below in more detail.

When dealing with fuzzy sets according to CFL, the universal set $\mathbb U$ (the set to which each element susceptible to being considered belongs, as will be accepted, with a weight of truth of 1) must be specified first. The membership of an element $x_i$, for $i = 1, 2, 3, \dots$, in the set $\mathbb U$ will be specified by the proposition $q_{i, \mathbb U}$, and as indicated $w(q_{i, \mathbb U}) = 1$.

According to CFL, the empty set $\varnothing$ is the set to which each element susceptible to being considered within the frame of the given $\mathbb U$ also belongs, but the weight of truth of each proposition that establishes the membership of each of those elements in $\varnothing$ is equal to 0. The proposition $x_i \in \varnothing$ will be symbolized by $q_{i, \varnothing}$, for $i = 1, 2, 3, \dots$, and the weight of truth of that proposition will be symbolized by $w(q_{i, \varnothing})$. According to that, $w(q_{i, \varnothing}) = 0$.

As seen, the assignment of weights of truth to the propositions that establish the membership of elements in the sets  
$\mathbb U$ and $\varnothing$ does not require the intervention of experts in the field covered.

Each of the elements belonging to $\mathbb U$ also belong to each subset $C_j$, for $j = 1, 2, 3, \dots$, of $\mathbb U$; but for that $C_j$ to be different from $\mathbb U$ at least one of the $w(q_{i, j})$ must be different from 1, and for that $C_j$ to be different from $\varnothing$ at least one of the $w(q_{i, j})$ must be different from 0.

Consider any set $C_j = 1, 2, 3, \dots$, characterized within the frame of a certain $\mathbb U$. The complement of that set is symbolized as $\xvec{C}_j$.

Just as $q_{i, j}$ is the abbreviated form of the proposition $x_i \in C_j$, $q_{i, \xvec{j}}$ is the abbreviated form of the proposition $x_i \in \xvec{C}_j$, 

Given $w(q_{i, j})$ (i.e., the weight of truth of $(q_{i, j})$), $w(q_{i, \xvec{j}})$ (i.e., the weight of truth of $q_{i, \xvec{j}}$) can be computed: $w(q_{i, \xvec{j}}) = w(\overline{q}_{i, \xvec{j}}) = 1 - w(q_{i, j})$.

Thus, for example, if $w(q_{i, j}) = 0.47$, then $w(q_{i, \xvec{j}}) = 1- 0.47 = 0.53$.

The complement of $\mathbb U$ will be symbolized as $\xvec{\mathbb U}$. The abbreviated form of $x_i \in \xvec{\mathbb U}$ is $q_{i, \xvec{\mathbb U}}\;$. Given that $w(q_{i, \mathbb U}) = 1$, $w(q_{i, \xvec{\mathbb U}}) = 1- w(q_{i,\mathbb U}) = 1-1 = 0$. That is, for any $x$ belonging to $\xvec{\mathbb U}$, the equality $w(q_{i, \xvec{\mathbb U}}) = 0$ is valid. Note that, as in the theory of classical sets, $\xvec{\varnothing} = 0$ also in fuzzy set theory, according to CFL.

The proposition $x_i \in \varnothing$ ia abbreviated as $q_{i, \varnothing}$. The complement of $\varnothing$ is symbolized as $\xvec{\varnothing}$. The proposition $x_i \in \xvec{\varnothing}$ is abbreviated as $q_{i, \xvec\varnothing}$. Given that $w(q_{i, \varnothing}) = 0, w(q_{i, \varnothing}) = w(\overline{q}_{i, \varnothing}) = 1 - w(q_{i, \varnothing}) = 1 - 0 = 1$. Note that in the theory of classical sets, $\xvec{\varnothing} = \mathbb U$. If in the above equation the operation of complementation of the equality is carried out, one obtains $\xvec{\xvec{\varnothing}} = \xvec{\mathbb U} = \varnothing$. That is, the double complementation of $\varnothing$ is equal to $\varnothing$. If the operation of complementation of both members of the equality $\xvec{\mathbb U} = \varnothing$ is carried out, one obtains $\xvec{\xvec{\mathbb U}} = \xvec{\varnothing} = \mathbb U$. That is, the double complementation of $\mathbb U$ is equal to $\mathbb U$.

In general, for any set $C_j$, for $j = 1, 2, 3, \dots$, in CFL, just as in the theory of classical sets, the following equality is valid: $\xvec{\xvec{C}}_j = C_j$.

In \cite{a1} consideration was given to some operations with fuzzy sets according to CFL. The coverage of those operations will be expanded and presented in greater detail below.

Consider the universal set $\mathbb U$ to which 5 elements belong: $x_1, x_2, x_3 , x_4$ and $x_5$. Those elements are people who will be considered from 2 points of view: 1) as chess players, and 2) as wealthy individuals.

It will be accepted that $x_2$ is an outstanding professional chess player, a grandmaster with possibilities of becoming the next world champion of that game/science.

The element $x_5$ is a chess fan whose level, according to chess experts, is that of a category 3 player.

The elements $x_1, x_3 $ and $x_4$ are individuals who are entirely unfamiliar with chess rules.

Note that a subset of the $\mathbb U$ considered is that of chess players in that  $\mathbb U$. Admit that according to the criterion of chess experts the weights of truth of $q_{1, 1}, q_{2, 1}, q_{3, 1}, q_{4, 1}$ and $q_{5, 1}$, are respectively, 
$w(q_{1, 1}) = 0, w(q_{2, 1}) = 1, w(q_{3, 1}) = 0, w(q_{4, 1})$ = 0, and $w(q_{5, 1}) = 0.4$.

The subset $C_2$ of the $\mathbb U$ considered is that of wealthy individuals. Admit that, according to the criterion of economists and experts in finance, the weights of truth of $q_{1, 1}, q_{2, 1}, q_{3, 1}, q_{4, 1}$, and $q_{5, 1}$, are respectively, $w(q_{1, 2}) = 0.9, w(q_{2, 2}) = 0.8, w(q_{3, 2}) = 0.7, w(q_{4, 2})$ = 0, and $w(q_{5, 2}) = 0.6$.

Attention will be given below to some operations carried out with the sets $\mathbb U, C_1$ and $C_2$ just addressed.

First, those sets and their complementary sets will be presented as column vectors. Second, the results obtained from each of the operations carried out will be presented. To obtain those results, operations should be carried out on the propositions corresponding to the different elements belonging to those sets. Anyone can carry out those operations, if desired. Third, as an example of these operations, detailed information will be provided after each result obtained from operations on fuzzy sets, according to CFL, only for the operations corresponding to those propositions concerning the element $x_5$.

\begin{align*}
\mathbb U = \begin{array}{|cc|}
x_1;  & w(q_{1, \mathbb U}) = 1 \\
x_2;  & w(q_{2, \mathbb U}) = 1 \\
x_3;  & w(q_{3, \mathbb U}) = 1 \\
x_4;  & w(q_{4, \mathbb U}) = 1 \\
x_5;  & w(q_{5, \mathbb U}) = 1 \\
\end{array} \hspace{.1in} ; \varnothing = \xvec{\mathbb U} = \hspace{.1in}  \begin{array}{|cc|}
x_1;  & w(q_{1,\varnothing}) = 0  \\
x_2;  & w(q_{2,\varnothing}) = 0 \\
x_3;  & w(q_{3,\varnothing}) = 0 \\
x_4;  & w(q_{4,\varnothing}) = 0 \\
x_5;  & w(q_{5,\varnothing}) = 0 \\
\end{array}\end{align*}

\noindent $ q_{5,\mathbb U}: x_5 \in \mathbb U; w(q_{5,\mathbb U}) = 1$\\
\noindent $ q_{5,\varnothing}: x_5 \in \varnothing; w(q_{5, \varnothing}) = w(q_{5,\xvec{\mathbb U}}) = 1 -  w(q_{5,\mathbb U}) = 1- 1 = 0$

\begin{align*}
C_1 = \begin{array}{|cccl|}
x_1;  & w(q_{1, 1}) &=& 0 \\
x_2;  & w(q_{2, 1}) &=& 1 \\
x_3;  & w(q_{3, 1}) &=& 0 \\
x_4;  & w(q_{4, 1}) &=& 0 \\
x_5;  & w(q_{5, 1}) &=& 0.4 \\
\end{array} \hspace{.1in} ; \xvec{C}_1 = \hspace{.1in}  \begin{array}{|cccl|}
x_1;  & w(q_{1,\xvec{1}}) &=& 1  \\
x_2;  & w(q_{2,\xvec{1}}) &=& 0 \\
x_3;  & w(q_{3,\xvec{1}}) &=& 1 \\
x_4;  & w(q_{4,\xvec{1}}) &=& 1 \\
x_5;  & w(q_{5,\xvec{1}}) &=& 0.6 \\
\end{array}\end{align*}

\noindent $ q_{5,1}: x_5 \in C_2; w(q_{5,1}) = 0.4$\\
\noindent $ q_{5,\xvec{1}}: x_5 \in \xvec{C}_1; w(q_{5,\xvec{1}}) = 1 -  w(q_{5,1}) = 1- 0.4 = 0.6$

\begin{align*}
C_2 = \begin{array}{|cccl|}
x_1;  & w(q_{1, 2}) &=& 0.9 \\
x_2;  & w(q_{2, 2}) &=& 0.8 \\
x_3;  & w(q_{3, 2}) &=& 0.7 \\
x_4;  & w(q_{4, 2}) &=& 0 \\
x_5;  & w(q_{5, 2}) &=& 0.6 \\
\end{array} \hspace{.1in} ; \xvec{C}_2 = \hspace{.1in}  \begin{array}{|cccl|}
x_1;  & w(q_{1,\xvec{2}}) &=& 0.1  \\
x_2;  & w(q_{2,\xvec{2}}) &=& 0.2 \\
x_3;  & w(q_{3,\xvec{2}}) &=& 0.3 \\
x_4;  & w(q_{4,\xvec{2}}) &=& 1 \\
x_5;  & w(q_{5,\xvec{2}}) &=& 0.4 \\
\end{array}\end{align*}

\noindent $ q_{5,2}: x_5 \in C_2; w(q_{5,2}) = 0.6$\\
\noindent $ q_{5,\xvec{2}}: x_5 \in \xvec{C}_2; w(q_{5,\xvec{2}}) = 1 -  w(q_{5,2}) = 1- 0.6 = 0.4$

\begin{align*}
C_1 \cup C_2 = \begin{array}{|cccl|}
x_1;  & w(q_{1,{1}}) &=& 0 \\
x_2;  & w(q_{2,{1}}) &=& 1 \\
x_3;  & w(q_{3, {1}}) &=& 0 \\
x_4;  & w(q_{4, {1}}) &=& 0 \\
x_5;  & w(q_{5, {1}}) &=& 0.4 \\
\end{array} \hspace{.1in} \cup \hspace{.1in}  \begin{array}{|cccl|}
x_1;  & w(q_{1,2}) &=& 0.9   \\
x_2;  & w(q_{2,2}) &=& 0.8 \\
x_3;  & w(q_{3,2}) &=& 0.7 \\
x_4;  & w(q_{4,2}) &=& 0 \\
x_5;  & w(q_{5,2}) &=& 0.6 \\
\end{array} = \\  \begin{array}{|cccl|}
x_1;  & w(q_{1,1} \lor q_{1,{2}}) &=& 0.9  \\
x_2;  & w(q_{2,1} \lor q_{2,{2}}) &=& 1 \\
x_3;  & w(q_{3,1} \lor q_{3,{2}}) &=& 0.7 \\
x_4;  & w(q_{4,1} \lor q_{4,{2}}) &=& 0 \\
x_5;  & w(q_{5,1} \lor q_{5,{2}}) &=& 0.76 \\
\end{array}\end{align*}

\noindent $ q_{5,1}: x_5 \in C_1; w(q_{5,1}) = 0.4$\\
$ q_{5,2}: x_5 \in C_2; w(q_{5,2}) = 0.6$\\

\begin{tabular}{c|c||c|l}
$q_{5,1}$ & $q_{5,2}$ & $q_{5,1} \lor q_{5,2}$ & $w(q_{5,1} \lor q_{5,2}) = \displaystyle\sum_{i=1}^{4}S_i$\\
\midrule
0 & 0 & 0 & $S_1 = 0 $\\
0 & 1 & 1 & $S_2 = w(\overline{q}_{5,1}) \cdot w(q_{5,2}) $\\   
1 & 0 & 1 & $S_3 = w(q_{5,1}) \cdot w(\overline{q}_{5,2}) $\\  
1 & 1 & 1 & $S_4 = w(q_{5,1}) \cdot w(q_{5,2}) $\\ 
\end{tabular}\\
\\

\noindent $w(q_{5,1} \lor q_{5,2}) = 0 + (0.6) \cdot (0.6) + (0.4) \cdot (0.4) + (0.6) \cdot (0.4) = 0.36 + 0.16 + 0.24 = 0.76$\\

\noindent That is, $w(q_5,1 \lor q_5,2) = 0.76$.

\begin{align*}
C_1 \cap C_2 = \begin{array}{|cccl|}
x_1;  & w(q_{1,{1}}) &=& 0   \\
x_2;  & w(q_{2,{1}}) &=&1 \\
x_3;  & w(q_{3, {1}}) &=& 0 \\
x_4;  & w(q_{4, {1}}) &=& 0 \\
x_5;  & w(q_{5, {1}}) &=& 0.4 \\
\end{array} \hspace{.1in}  \cap \hspace{.1in}  \begin{array}{|cccl|}
x_1;  & w(q_{1,2}) &=& 0.9   \\
x_2;  & w(q_{2,2}) &=& 0.8 \\
x_3;  & w(q_{3,2}) &=& 0.7 \\
x_4;  & w(q_{4,2}) &=& 0 \\
x_5;  & w(q_{5,2}) &=& 0.6 \\
\end{array} = \\  \begin{array}{|cccl|}
x_1;  & w(q_{1,1} \land q_{1,{2}}) &=& 0  \\
x_2;  & w(q_{2,1} \land q_{2,{2}}) &=& 0.8 \\
x_3;  & w(q_{3,1} \land q_{3,{2}}) &=& 0 \\
x_4;  & w(q_{4,1} \land q_{4,{2}}) &=& 0 \\
x_5;  & w(q_{5,1} \land q_{5,{2}}) &=& 0.24 \\
\end{array}\end{align*}

\noindent $ q_{5,1}: x_5 \in C_1; w(q_{5,1}) = 0.4$\\
$ q_{5,2}: x_5 \in C_2; w(q_{5,2}) = 0.6$\\

\begin{tabular}{c|c||c|l}
$q_{5,1}$ & $q_{5,2}$ & $q_{5,1} \land q_{5,2}$ & $w(q_{5,1} \land q_{5,2}) = \displaystyle\sum_{i=1}^{4}S_i$\\
\midrule
0 & 0 & 0 & $S_1 = 0 $\\
0 & 1 & 0 & $S_2 = 0 $\\   
1 & 0 & 0 & $S_3 = 0 $\\  
1 & 1 & 1 & $S_4 = w(q_{5,1}) \cdot w(q_{5,2}) = (0.4) \cdot (0.6) = 0.24$
\end{tabular}\\
\\
\\
\noindent $w(q_{5,1} \land q_{5,2}) = 0 + 0 + 0 + (0.4) \cdot (0.6) = 0.24$

\begin{align*}
C_1 \dot{\cup} C_2 = \begin{array}{|cccl|}
x_1;  & w(q_{1,{1}}) &=& 0   \\
x_2;  & w(q_{2,{1}}) &=& 1 \\
x_3;  & w(q_{3, {1}}) &=& 0 \\
x_4;  & w(q_{4, {1}}) &=& 0 \\
x_5;  & w(q_{5, {1}}) &=& 0.4 \\
\end{array} \hspace{.1in}  \dot{\cup} \hspace{.1in}  \begin{array}{|cccl|}
x_1;  & w(q_{1,2}) &=& 0.9   \\
x_2;  & w(q_{2,2}) &=& 0.8 \\
x_3;  & w(q_{3,2}) &=& 0.7 \\
x_4;  & w(q_{4,2}) &=& 0 \\
x_5;  & w(q_{5,2}) &=& 0.6 \\
\end{array} = \\  \begin{array}{|cccl|}
x_1;  & w(q_{1,1} \dot{\lor} q_{1,{2}}) &=& 0.9  \\
x_2;  & w(q_{2,1} \dot{\lor} q_{2,{2}}) &=& 0.2 \\
x_3;  & w(q_{3,1} \dot{\lor} q_{3,{2}}) &=& 0.7 \\
x_4;  & w(q_{4,1} \dot{\lor} q_{4,{2}}) &=& 0 \\
x_5;  & w(q_{5,1} \dot{\lor} q_{5,{2}}) &=& 0.52 \\
\end{array}\end{align*}

\noindent $ q_{5,1}: x_5 \in C_1; w(q_{5,1}) = 0.4$\\
$ q_{5,2}: x_5 \in C_2; w(q_{5,2}) = 0.6$\\

\begin{tabular}{c|c||c|l}
$q_{5,1}$ & $q_{5,2}$ & $q_{5,1} \dot\lor q_{5,2}$ & $w(q_{5,1} \dot\lor q_{5,2}) = \displaystyle\sum_{i=1}^{4}S_i$\\
\midrule
0 & 0 & 0 & $S_1 = 0 $\\
0 & 1 & 1 & $S_2 = w(\overline{q}_{5,1}) \cdot w(q_{5,2}) $\\   
1 & 0 & 1 & $S_3 = w(q_{5,1}) \cdot w(\overline{q}_{5,2}) $\\  
1 & 1 & 0 & $S_4 = 0 $\\ 
\end{tabular}\\
\\

\noindent $w(q_{5,1} \dot\lor q_{5,2}) = 0 + (0.6) \cdot (0.6) + (0.4) \cdot (0.4) + 0= 0.36 + 0.16 = 0.52$

\begin{align*}
C_1 \naturalto C_2 = \begin{array}{|cccl|}
x_1;  & w(q_{1,{1}}) &=& 0   \\
x_2;  & w(q_{2,{1}}) &=& 1 \\
x_3;  & w(q_{3, {1}}) &=& 0 \\
x_4;  & w(q_{4, {1}}) &=& 0 \\
x_5;  & w(q_{5, {1}}) &=& 0.4 \\
\end{array} \hspace{.1in} \naturalto \hspace{.1in}  \begin{array}{|cccl|}
x_1;  & w(q_{1,2}) &=& 0.9   \\
x_2;  & w(q_{2,2}) &=& 0.8 \\
x_3;  & w(q_{3,2}) &=& 0.7 \\
x_4;  & w(q_{4,2}) &=& 0 \\
x_5;  & w(q_{5,2}) &=& 0.6 \\
\end{array} = \\  \begin{array}{|cccl|}
x_1;  & w(q_{1,1} \to q_{1,{2}}) &=& 1 \\
x_2;  & w(q_{2,1} \to q_{2,{2}}) &=& 0.8 \\
x_3;  & w(q_{3,1} \to q_{3,{2}}) &=& 1 \\
x_4;  & w(q_{4,1} \to q_{4,{2}}) &=& 1 \\
x_5;  & w(q_{5,1} \to q_{5,{2}}) &=& 0.84 \\
\end{array}\end{align*}

\noindent $ q_{5,1}: x_5 \in C_1; w(q_{5,1}) = 0.4$\\
$ q_{5,2}: x_5 \in C_2; w(q_{5,2}) = 0.6$\\

\begin{tabular}{c|c||c|l}
$q_{5,1}$ & $q_{5,2}$ & $q_{5,1} \to q_{5,2}$ & $w(q_{5,1} \to q_{5,2}) = \displaystyle\sum_{i=1}^{4}S_i$\\
\midrule
0 & 0 & 1 & $S_1 = w(\overline{q}_{5,1}) \cdot w(\overline{q}_{5,2}) $\\
0 & 1 & 1 & $S_2 = w(\overline{q}_{5,1}) \cdot w(q_{5,2}) $\\   
1 & 0 & 0 & $S_3 = 0$\\
1 & 1 & 1 & $S_4 = w(q_{5,1}) \cdot w(q_{5,2}) $\\
\end{tabular}\\
\\
\\
\noindent $w(q_{5,1} \to q_{5,2}) = (0.6) \cdot (0.4) + (0.6) \cdot (0.6) + 0 + (0.4) \cdot (0.6) = 0.24 + 0.36 + 0.24 = 0.84$
\begin{align*}
C_2 \naturalto C_1 = \begin{array}{|cccl|}
x_1;  & w(q_{1,{2}}) &=& 0.9  \\
x_2;  & w(q_{2,{2}}) &=& 0.8 \\
x_3;  & w(q_{3, {2}}) &=& 0.7 \\
x_4;  & w(q_{4, {2}}) &=& 0 \\
x_5;  & w(q_{5, {2}}) &=& 0.6 \\
\end{array} \hspace{.1in} \naturalto \hspace{.1in}  \begin{array}{|cccl|}
x_1;  & w(q_{1,1}) &=& 0   \\
x_2;  & w(q_{2,1}) &=& 1 \\
x_3;  & w(q_{3,1}) &=& 0 \\
x_4;  & w(q_{4,1}) &=& 0 \\
x_5;  & w(q_{5,1}) &=& 0.4 \\
\end{array} = \\  \begin{array}{|cccl|}
x_1;  & w(q_{1,2} \to q_{1,{1}}) &=& 0.1 \\
x_2;  & w(q_{2,2} \to q_{2,{1}}) &=& 1 \\
x_3;  & w(q_{3,2} \to q_{3,{1}}) &=& 0.3 \\
x_4;  & w(q_{4,2} \to q_{4,{1}}) &=& 1 \\
x_5;  & w(q_{5,2} \to q_{5,{1}}) &=& 0.64 \\
\end{array}\end{align*}

\noindent $ q_{5,2}: x_5 \in C_2; w(q_{5,1}) = 0.6$\\
$ q_{5,1}: x_5 \in C_1; w(q_{5,2}) = 0.4$\\

\begin{tabular}{c|c||c|l}
$q_{5,2}$ & $q_{5,1}$ & $q_{5,2} \to q_{5,1}$ & $w(q_{5,2} \to q_{5,1}) = \displaystyle\sum_{i=1}^{4}S_i$\\
\midrule
0 & 0 & 1 & $S_1= w(\overline{q}_{5,2}) \cdot w(\overline{q}_{5,1}) $\\
0 & 1 & 1 & $S_2= w(\overline{q}_{5,2}) \cdot w(q_{5,1}) $\\   
1 & 0 & 0 & $S_3= 0$\\
1 & 1 & 1 & $S_4= w(q_{5,2}) \cdot w(q_{5,1}) $\\
\end{tabular}\\
\\

\noindent $w(q_{5,2} \to q_{5,1}) = (0.4) \cdot (0.6) + (0.4) \cdot (0.4) + 0 + (0.6) \cdot (0.4) = 0.24 + 0.16 + 0.24 = 0.64$

\begin{align*}
C_1 \naturaltolr C_2 = \begin{array}{|cccl|}
x_1;  & w(q_{1,{1}}) &=& 0  \\
x_2;  & w(q_{2,{1}}) &=& 1 \\
x_3;  & w(q_{3, {1}}) &=& 0 \\
x_4;  & w(q_{4, {1}}) &=& 0 \\
x_5;  & w(q_{5, {1}}) &=& 0.4 \\
\end{array} \hspace{.1in} \naturaltolr \hspace{.1in}  \begin{array}{|cccl|}
x_1;  & w(q_{1,2}) &=& 0.9   \\
x_2;  & w(q_{2,2}) &=& 0.8 \\
x_3;  & w(q_{3,2}) &=& 0.7 \\
x_4;  & w(q_{4,2}) &=& 0 \\
x_5;  & w(q_{5,2}) &=& 0.6 \\
\end{array} = \\  \begin{array}{|cccl|}
x_1;  & w(q_{1,1} \longleftrightarrow q_{1, 2}) &=& 0.1 \\
x_2;  & w(q_{2,1} \longleftrightarrow q_{2, 2}) &=& 0.8 \\
x_3;  & w(q_{3,1} \longleftrightarrow q_{3, 2}) &=& 0.3 \\
x_4;  & w(q_{4,1} \longleftrightarrow q_{4, 2}) &=& 1 \\
x_5;  & w(q_{5,1} \longleftrightarrow q_{5, 2}) &=& 0.48 \\
\end{array}\end{align*}

\noindent $ q_{5,1}: x_5 \in C_1; w(q_{5,1}) = 0.4$\\
$ q_{5,2}: x_5 \in C_2; w(q_{5,2}) = 0.6$\\

\begin{tabular}{c|c||c|l}
$q_{5,1}$ & $q_{5,2}$ & $q_{5,1} \longleftrightarrow q_{5,2}$ & $w(q_{5,1} \longleftrightarrow q_{5,2}) = \displaystyle\sum_{i=1}^{4}S_i$\\
\midrule
0 & 0 & 1 & $S_1 = w(\overline{q}_{5,2}) \cdot w(\overline{q}_{5,1}) $\\
0 & 1 & 0 & $S_2 = 0$\\
1 & 0 & 0 & $S_3 = 0$\\
1 & 1 & 1 & $S_4 = w(q_{5,1}) \cdot w(q_{5,2}) $\\
\end{tabular}\\
\\

\noindent $w(q_{5,1} \longleftrightarrow q_{5,2}) = (0.6) \cdot (0.4) + 0 + 0+ (0.4) \cdot (0.6) = 0.24 + 0.24 = 0.48$\\
\begin{align*}
C_1 \uparrowbarred C_2 = \begin{array}{|cccl|}
x_1;  & w(q_{1,{1}}) &=& 0  \\
x_2;  & w(q_{2,{1}}) &=& 1 \\
x_3;  & w(q_{3, {1}}) &=& 0 \\
x_4;  & w(q_{4, {1}}) &=& 0 \\
x_5;  & w(q_{5, {1}}) &=& 0.4 \\
\end{array} \hspace{.1in} \uparrowbarred \hspace{.1in}  \begin{array}{|cccl|}
x_1;  & w(q_{1,2}) &=& 0.9   \\
x_2;  & w(q_{2,2}) &=& 0.8 \\
x_3;  & w(q_{3,2}) &=& 0.7 \\
x_4;  & w(q_{4,2}) &=& 0 \\
x_5;  & w(q_{5,2}) &=& 0.6 \\
\end{array} = \\  \begin{array}{|cccl|}
x_1;  & w(q_{1,1} \uparrow q_{1,{2}}) &=& 1 \\
x_2;  & w(q_{2,1} \uparrow q_{2,{2}}) &=& 0.2 \\
x_3;  & w(q_{3,1} \uparrow q_{3,{2}}) &=& 1 \\
x_4;  & w(q_{4,1} \uparrow q_{4,{2}}) &=& 1 \\
x_5;  & w(q_{5,1} \uparrow q_{5,{2}}) &=& 0.76 \\
\end{array}\end{align*}

\noindent $ q_{5,1}: x_5 \in C_1; w(q_{5,1}) = 0.4$\\
$ q_{5,2}: x_5 \in C_2; w(q_{5,2}) = 0.6$\\

\begin{tabular}{c|c||c|c}
$q_{5,1}$ & $q_{5,2}$ & $q_{5,1} \uparrow q_{5,2}$ & $w(q_{5,1} \uparrow q_{5,2}) = \displaystyle\sum_{i=1}^{4}S_i$\\
\midrule
0 & 0 & 1 & $S_1 = w(\overline{q}_{5,1}) \cdot w(\overline{q}_{5,2}) $\\
0 & 1 & 1 & $S_2 = w(\overline{q}_{5,1}) \cdot w(q_{5,2}) $\\   
1 & 0 & 1 & $S_3 = w(q_{5,1}) \cdot w(\overline{q}_{5,2}) $\\
1 & 1 & 0 & $S_4 = 0 $\\
\end{tabular}\\
\\

\noindent $w(q_{5,1} \uparrow q_{5,2}) = (0.6) \cdot (0.4) + (0.6) \cdot (0.6) + (0.4) \cdot (0.4) + 0 = 0.24 + 0.36 + 0.16 = 0.76$
\begin{align*}
C_1 \downarrowbarred C_2 = \begin{array}{|cccl|}
x_1;  & w(q_{1,{1}}) &=  &0  \\
x_2;  & w(q_{2,{1}}) &=& 1 \\
x_3;  & w(q_{3, {1}}) &= &0 \\
x_4;  & w(q_{4, {1}}) &= &0 \\
x_5;  & w(q_{5, {1}}) &=& 0.4 \\
\end{array} \hspace{.1in} \downarrowbarred \hspace{.1in}  \begin{array}{|cccl|}
x_1;  & w(q_{1,2}) &= &0.9   \\
x_2;  & w(q_{2,2}) &= &0.8 \\
x_3;  & w(q_{3,2}) &= &0.7 \\
x_4;  & w(q_{4,2}) &= &0 \\
x_5;  & w(q_{5,2}) &=& 0.6 \\
\end{array} = \\  \begin{array}{|cccl|}
x_1;  & w(q_{1,1} \downarrow q_{1,{2}}) &=& 0.1 \\
x_2;  & w(q_{2,1} \downarrow q_{2,{2}}) &= &0 \\
x_3;  & w(q_{3,1} \downarrow q_{3,{2}}) &=& 0.3 \\
x_4;  & w(q_{4,1} \downarrow q_{4,{2}}) &= &1 \\
x_5;  & w(q_{5,1} \downarrow q_{5,{2}}) &= &0.24 \\
\end{array}\end{align*}

\begin{tabular}{c|c||c|l}
$q_{5,1}$ & $q_{5,2}$ & $q_{5,1} \downarrow q_{5,1}$ & $w(q_{5,2} \downarrow q_{5,2}) = \displaystyle\sum_{i=1}^{4}S_i$\\
\midrule
0 & 0 & 1 & $S_1 = w(\overline{q}_{5,1}) \cdot w(\overline{q}_{5,2}) $\\
0 & 1 & 0 & $S_2 = 0$\\
1 & 0 & 0 & $S_3 = 0$\\
1 & 1 & 0 & $S_4 = 0$\\
\end{tabular}\\
\\

\noindent $w(q_{5,1} \downarrow q_{5,2}) = (0.6) \cdot (0.4) + 0 + 0+ 0 = 0.24 $\\

\subsection{Fuzzy Sets, according to CFL, to Which Infinite Elements Belong}

Suppose that infinite elements belong to a fuzzy set $C$ according to CFL. Consider, for example, that each of those elements can be identified as one of the numbers belonging to the interval of the real variable $x$ between $x = 0$ and $x = 10$. Given that in that interval there are infinite elements, to assign a weight of truth  to each proposition that affirms the membership of one of those variables $x$ in that interval, the value of truth can be established as a function of $x$ and of $C$: $w(x \in C); 0 \leq x \leq 10$. The proposition $(x \in C)$ will be abbreviated as $q_{x, C}$. Thus, the weight of truth of the proposition $x \in C$ can be expressed as $w(q_{x, C}); 0 \leq x \leq 10$.

If it is accepted that for any of these elements $x$ the corresponding weight of truth equals $1$, then $C$ will be equal to $\mathbb U$; that is, equal to the universal set considered: $C = \mathbb U$. According to the above, $w(q_{x, \mathbb U}) = 1; 0 \leq x \leq 10$.

The complementary set of $\mathbb U$ is the corresponding empty set: $\xvec{\mathbb U} = \varnothing$. Each element $x$ belonging to $\mathbb U$ also belongs to $\varnothing$, and the proposition that, when abbreviated, expresses the membership of that $x$ in this latter set is $(x \in \varnothing)$. The weight of truth of this proposition for any $x$ belonging to $\mathbb U$ is equal to $0$: $w(q_{x, \varnothing}) = 0$; $0 \leq x \leq 10$. In effect, with $\varnothing = \xvec{\mathbb U}$, the result is $w(q_{x, \varnothing}) = w(\xvec{q}_{x, \mathbb U}) = 1 - w(q_{x, \mathbb U}) = 1 - 1 = 0$. 

\begin{figure}[H]
\begin{tikzpicture}
\draw[thin,->] (0,0) -- (10.5,0) node[anchor=north west] { $x$ };
\draw[thin,->] (0,0) -- (0,4.5) node[anchor=south east] {  };
\draw (1pt,4) -- (-1pt,4) node[anchor=east] { 1 };
\draw(0,1pt) -- (0,-1pt) node[anchor=north] {0};
\draw (5,1pt) -- (5,-1pt) node[anchor=north] {5 };
\draw (10,1pt) -- (10,-1pt) node[anchor=north] { 10 };
\draw[very thick] (0,4) -- (10,4);
\draw[->] (8,3.9) to[bend right] (9,3) node[anchor=west] {  $w(q_{x,\mathbb{U}}) =1$ };
\draw[very thick] (0,0) -- (10,0);
\draw[->] (8,0.1) to[bend left] (9,1)  node[anchor=west] {  $w(q_{x,\varnothing}) =0$ } ;
\end{tikzpicture}
\caption{Weights of truth $w(q_{x, \mathbb U})$ and $w(q_{x, \varnothing})$ in the interval considered $0 \leq x \leq 10$}
\label{f28}
\end{figure}
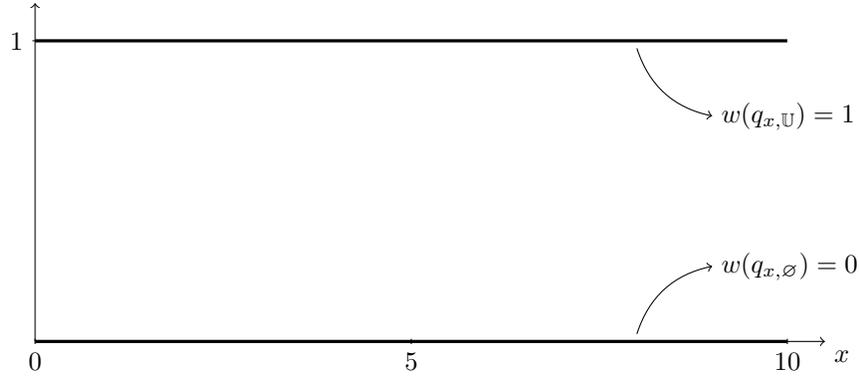

The weights of truth $w(q_{x, \mathbb U})$ and $w(q_{x, \varnothing})$ in the interval considered $0 \leq x \leq 10$ are shown in figure \ref{f28}.

Two sets denominated $C_1$ and $C_2$ will be characterized below within the frame of the $\mathbb U$ considered.\\

\noindent $q_{x, 1}: x  \in C_1$

If $(0 \leq x \leq 5)$, then $w(q_{x, 1}) = \frac {1}{5} \cdot x$.

If $(5 \leq x \leq 10)$, then $w(q_{x, 1}) = 1 - \frac {1}{5}(x-5)$.\\

\noindent$(q_{x, 2}): x \in C_2$

If $(0 \leq x < 3)$, then $w(q_{x, 2}) = 0$.

If $(3 \leq x < 7)$, then $w(q_{x, 2}) = 0.8$.

If $(7 \leq x \leq 10)$, then $w(q_{x, 2}) = 0$.\\

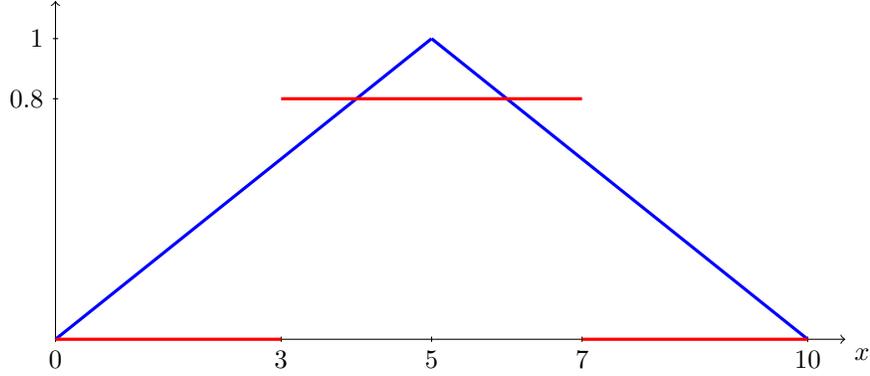
\begin{figure}[H]
\begin{tikzpicture}
\draw[thin,->] (0,0) -- (10.5,0) node[anchor=north west] { $x$ };
\draw[thin,->] (0,0) -- (0,4.5) node[anchor=south east] {  };
\draw (1pt,4) -- (-1pt,4) node[anchor=east] { 1 };
\draw (1pt, 3.2) -- (-1pt, 3.2) node[anchor=east] {0.8};
\draw(0,1pt) -- (0,-1pt) node[anchor=north] {0};
\draw (5,1pt) -- (5,-1pt) node[anchor=north] {5 };
\draw (10,1pt) -- (10,-1pt) node[anchor=north] { 10 };

\draw (3,1pt) -- (3,-1pt) node[anchor=north] {3 };
\draw (7,1pt) -- (7,-1pt) node[anchor=north] { 7 };

\draw[very thick,blue] (0,0) -- (5,4);
\draw[very thick,blue] (5,4) -- (10,0);
\draw[very thick,red] (0,0) -- (3,0);
\draw[very thick,red] (3,3.2) -- (7,3.2);
\draw[very thick,red] (7,0) -- (10,0);
\end{tikzpicture}
\caption{Weights of truth in the interval considered $0 \leq x \leq 10$: $w(q_{x, 1})$ -- blue; and $w(q_{x, 2})$ -- red}
\label{f29}
\end{figure}

The weights of truth $w(q_{x, 1})$ and $w(q_{x, 2})$ in the interval $0 \leq x \leq 10$ are characterized in figure \ref{f29}.

The sets $C_3 = C_1 \cup C_2$ and $C_4 = C_1 \cap C_2$ will be considered below.

The weight of truth of the proposition $x \in C_3$ will be symbolized as $w(q_{x, 3})$.  The weight of truth of the proposition $x \in C_4$ will be symbolized as $w(q_{x, 4})$.

First, $w(q_{x, 3}) = w(q_{x, 1} \lor q_{x, 2})$ will be computed for subintervals of the interval $0 \leq x \leq 10$.

\begin{tabular}{c|c||c|l}
$q_{x, 1}$ & $q_{x, 2}$ & $q_{x, 1} \lor q_{x, 2}$ & $w(q_{2, 3}) = w(q_{x, 1} \lor q_{x, 2}) =\displaystyle\sum_{i=1}^{4}S_i$\\
\midrule
0 & 0 & 0 & $S_1 = 0 $\\
0 & 1 & 1 & $S_2 = w(\overline{q}_{x, 1}) \cdot w(q_{x, 2}) $\\   
1 & 0 & 1 & $S_3 = w(q_{x, 1}) \cdot w(\overline{q}_{x, 2}) $\\
1 & 1 & 1 & $S_4 = w(q_{x, 1}) \cdot w(q_{x, 2}) $\\
\end{tabular}\\
\\

If $(0 \leq x < 3)$, then $w(q_{x, 1}) = \frac {1}{5}x$ and $w(q_{x, 2}) = 0$. Therefore, $w(q_{x, 3}) = 0 + (1 - \frac {1}{5}x) \cdot 0 + (\frac {1}{5}x) (1 - 0) + (\frac {1}{5}x) \cdot 0 = \frac {1}{5}x$. That is, if $(0 \leq x < 3)$, then $w(q_{x, 3}) = \frac {1}{5}x$.

If $(3 \leq x \leq 5)$, then $w(q_{x, 1}) = \frac {1}{5}x$ and $w(q_{x, 2}) = 0.8$. Therefore, $w(q_{x, 3}) = 0 + (1 - \frac{1}{5}x) \cdot (0.8) +(\frac{1}{5}x) \cdot (0.2) + (\frac{1}{5}x) \cdot (0.8) = 0.8 - (\frac{0.8}{5}) \cdot x + (\frac{0.2}{5}) \cdot x +(\frac{0.8}{5}) \cdot x  = 0.8 + (\frac{0.2}{5}) \cdot x  = 0.8 + (0.04) \cdot x$. That is, if $(3 \leq x \leq 5)$, then $w(q_{x, 3}) =  0.8 + (0.04) \cdot x$.

If $(5 < x \leq 7)$, then $w(q_{x, 1}) = 1 - \frac{1}{5}(x-5)$ and $w(q_{x, 2}) = 0.8$. Therefore, $w(q_{x, 3}) = 0 + \frac{1}{5} (x - 5) \cdot (0.8) + (1 - \frac{1}{5} (x - 5)) \cdot (0.2) + (1 - \frac{1}{5} (x - 5)) \cdot (0.8) = 1 - (\frac{0.2}{5}) \cdot (x - 5) = 1 - (0.04) \cdot (x - 5)$. That is, if $(5  < x \leq 7)$, then $w(q_{x, 3}) = 1 - (0.04) \cdot (x - 5)$.

If $(7 < x \leq 10)$, then $w(q_{x, 1}) = 1 - (\frac{1}{5}) \cdot (x - 5)$ and $w(q_{x, 2}) = 0$. Therefore, $w(q_{x, 3}) = 0 + ((\frac{1}{5} (x - 5)) \cdot 0 + (1 - \frac{1}{5} (x - 5)) \cdot 1 + (1 - \frac{1}{5} (x - 5)) \cdot 0 = 1 - \frac{1}{5}(x - 5) = 1 - (0.02) \cdot (x - 5)$. That is, if $7 < x \leq 10)$, then $w(q_{x, 3}) = 1 - (0.2) \cdot (x - 5)$.

Second, $w(q_{x, 4}) = w(q_{x, 1} \land q_{x, 2})$ will be computed for subintervals of the interval $0 \leq x \leq 10$.

\begin{tabular}{c|c||c|l}
$q_{x, 1}$ & $q_{x, 2}$ & $q_{x, 1} \land q_{x, 2}$ & $w(q_{x, 4}) = w(q_{x, 1} \land q_{x, 2}) =\displaystyle\sum_{i=1}^{4}S_i$\\
\midrule
0 & 0 & 0 & $S_1 = 0 $\\
0 & 1 & 0 & $S_2 = 0 $\\   
1 & 0 & 0 & $S_3 = 0$\\
1 & 1 & 1 & $S_4 = w(q_{x, 1}) \cdot w(q_{x, 2}) $\\
\end{tabular}\\
\\
\\
\indent If $(0 \leq x < 3)$, then $w(q_{x, 1}) = \frac {1}{5}x$ and $w(q_{x, 2}) = 0$. Therefore, $w(q_{x, 4}) = 0 + 0 + 0 + (\frac {1}{5}x) \cdot 0 = 0$. That is, if $(0 \leq x < 3)$, then $w(q_{x, 4}) = 0$.

If $(3 \leq x \leq 5)$, then $w(q_{x, 1}) = \frac {1}{5}x$ and $w(q_{x, 2}) = 0.8$. Therefore, $w(q_{x, 4}) = 0 + 0 + 0 + (\frac{1}{5}x) \cdot (0.8) = (0.16) \cdot x$. That is, if $(3 \leq x \leq 5)$, then $w(q_{x, 4}) =  (0.16) \cdot x$.

If $(5 < x \leq 7)$, then $w(q_{x, 1}) = 1 - \frac{1}{5}(x - 5)$ and $w(q_{x, 2}) = 0.8$. Therefore, if $(5 < x \leq 7)$, then $w(q_{x, 4}) = 0 + 0 + 0 + (1 - \frac{1}{5}(x - 5)) \cdot (0.8) = 1.6 - 0.16x$. That is, if $(5  < x \leq 7)$, then $w(q_{x, 4}) = 1. 6 - (0.16)x$.

If $(7 < x \leq 10)$, then $w(q_{x, 1}) = 1 - \frac{1}{5}(x - 5)$ and $w(q_{x, 2}) = 0$. Therefore, $w(q_{x, 4}) = 0 + 0 + 0 + (1 - \frac{1}{5}(x - 5)) \cdot  0 = 0$. That is, if $7 < x \leq 10)$, then $w(q_{x, 4}) = 0$.

\begin{figure}[H]
\begin{tikzpicture}
\draw[thin,->] (0,0) -- (10.5,0) node[anchor=north west] { $x$ };
\draw[thin,->] (0,0) -- (0,4.5) node[anchor=south east] {  };
\draw(0,1pt) -- (0,-1pt) node[anchor=north] {0};
\draw (1pt,4) -- (-1pt,4) node[anchor=east] { 1 };
\draw (1pt, 3.2) -- (-1pt, 3.2) node[anchor=east] {0.8};
\draw (5,1pt) -- (5,-1pt) node[anchor=north] {5 };
\draw (10,1pt) -- (10,-1pt) node[anchor=north] { 10 };

\draw (3,1pt) -- (3,-1pt) node[anchor=north] {3 };
\draw (7,1pt) -- (7,-1pt) node[anchor=north] { 7 };

\draw[very thick, black!60!green] (0,0)--(3,2.4);
\draw[very thick, black!60!green] (3,3.68)--(5,4);
\draw[very thick, black!60!green] (5,4)--(7,3.68);
\draw[very thick, black!60!green] (7,2.4)--(10,0);

\draw[very thick, black!50!brown] (0,0)--(3,0);
\draw[very thick, black!50!brown] (3,1.92)--(5,3.2);
\draw[very thick, black!50!brown] (5,3.2)--(7,1.92);
\draw[very thick, black!50!brown] (7,0)--(10,0);
\draw[->] (6,3.82) to[bend left] (8,3) node[anchor=north ] {  $w(q_{x,1} \lor q_{x,2}) $ };
\draw[->] (6,2.56) to[bend right] (9,1) node[anchor=west] {  $w(q_{x,1} \land q_{x,2}) $ };
\end{tikzpicture}
\caption{Weights of truth in the interval $0 \leq x \leq 10$: $w(q_{x, 3}) = w(q_{x, 1} \lor q_{x, 2})$ -- green; and $w(q_{x, 4})= w(q_{x, 1} \land q_{x, 2})$ -- brown}.
\label{f30}
\end{figure}
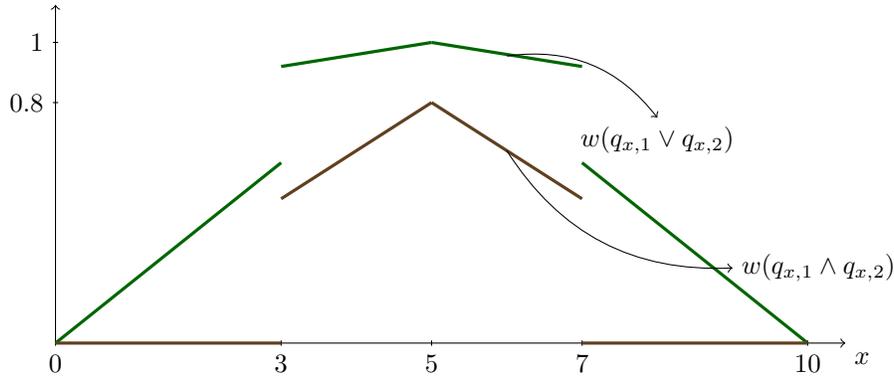

The weights of truth $w(q_{x, 3}) = w(q_{x, 1} \lor q_{x, 2})$ and $w(q_{x, 4}) = w(q_{x, 1} \land q_{x, 2})$ in the interval $0 \leq x \leq 10$ are shown in figure \ref{f30}.

It would be feasible to provide computations of the weights of truth of propositions corresponding to the sets resulting from other operations discussed above, such as the sets $C_5 = C_1 \dot{\cup} C_2, C_6, C_1 \naturalto C_2$, etc. It also would be possible to provide the results of operations with sets such that the propositions corresponding to their elements would have weights of truth that are functions different from those considered above, of those elements and of those sets to which they belong. Before presenting results of that type, it would be useful to provide software to facilitate those operations. This is addressed again in section 11.

\section{The Laws or Tautologies Valid for Classical Sets are also Valid for Both Finite and Infinite Fuzzy Sets according to CFL}

Consider any universal set $\mathbb U$ to which a finite or infinite number of elements belongs. Admit that one selects any element to be denominated $x$ belonging to that $\mathbb U$. To be selected, the only condition imposed on $x$ is the following: $x \in \mathbb U$.

Admit that $C_1$ and $C_2$ are two subsets whatsoever of that $\mathbb U$. The propositions that specify membership of $x$ in $C_1$ and $C_2$ are, respectively, the following: $x \in C_1$ abbreviated as $q_{x, 1}$, and $x \in C_2$ abbreviated as $q_{x, 2}$. Each of these propositions will have a certain weight of truth.

Consider, for example, the application of one of the De Morgan's laws to these propositions: \noindent $(\overline{q_{x, 1} \land q_{x, 2}} ) \longleftrightarrow (\overline{q}_{x, 1} \lor \overline {q}_{x, 2})$.

As seen above, $w((\overline {q_{x, 1} \land q_{x, 2}} ) \longleftrightarrow (\overline {q}_{x, 1} \lor \overline {q}_{x, 2})) = 1$, given that the weight of truth of any law --  or tautology -- in propositional calculus according to CFL is equal to 1.

Because $x$, according to the above, can be any element belonging to the $\mathbb U$ considered, the weight of truth of the proposition that affirms that each element of $\mathbb U$ belongs to the set $((\xvec{C_1 \cap C_2} ) \naturaltolr (\xvec{C_1} \cup \xvec{C_2} ))$  -- the one corresponding to that De Morgan's law -- is equal to 1. Therefore, that set is equal to $\mathbb U$; $((\xvec{C_1 \cap C_2} ) \naturaltolr (\xvec{C_1} \cup \xvec{C_2} ))$ = $\mathbb U$, an equality that constitutes the law of fuzzy set theory according to CFL corresponding to the De Morgan's law considered.

If consideration is given to any law $l(q_{x, 1}), q_{x, 2}, q_{x, 3}, \dots, q_{x, n})$ of propositional calculus in CFL with $n$ subsets of $\mathbb U$, one obtains -- as seen in section 6 -- the following result: $w(l(q_{x, 1}), q_{x, 2}, q_{x, 3}, \dots, q_{x, n})) = 1$. The set corresponding to that law is equal to $\mathbb U$.

\section{Classical Sets according to CFL}

Consider, within the framework of the classical theory of sets, the following universal set $\mathbb U$ to which 7 elements, denominated $x_1, x_2, x_3, x_4, x_5, x_6$ and $x_7$, belong: \noindent $\mathbb U = \{x_1, x_2, x_3, x_4, x_5, x_6$, $x_7 \}$.

Within the CFL framework, the following set corresponds to the above set:
\begin{align*}
\mathbb U = \begin{array}{|cc|}
x_1;  & w(q_{1, \mathbb U}) = 1 \\
x_2;  & w(q_{2, \mathbb U}) = 1 \\
x_3;  & w(q_{3, \mathbb U}) = 1 \\
x_4;  & w(q_{4, \mathbb U}) = 1 \\
x_5;  & w(q_{5, \mathbb U}) = 1 \\
x_6;  & w(q_{6, \mathbb U}) = 1 \\
x_7;  & w(q_{7, \mathbb U}) = 1 \\
\end{array} \end{align*}

In the theory of classical sets no element belongs to the empty set corresponding to the classical $\mathbb U$ considered. As explained, according to CFL each element belonging to $\mathbb U$ such that each $w(q_{i, \mathbb U}) = 0$, for $1, 2, 3, 4, 5, 6, 7$, belongs to that set. That is, the weight of truth of each proposition $q_{i, \varnothing}$ is equal to 0:
\begin{align*}
\varnothing = \begin{array}{|cc|}
x_1;  & w(q_{1, \varnothing}) = 0 \\
x_2;  & w(q_{2, \varnothing}) = 0 \\
x_3;  & w(q_{3, \varnothing}) = 0 \\
x_4;  & w(q_{4, \varnothing}) = 0 \\
x_5;  & w(q_{5, \varnothing}) = 0 \\
x_6;  & w(q_{6, \varnothing}) = 0 \\
x_7;  & w(q_{7, \varnothing}) = 0 \\
\end{array} \end{align*}

Consider the following subsets of $\mathbb U$ within the framework of the classical theory of sets:  $ C_1 = \{ 2, 3, 5, 7 \}; C_2 = \{1, 2, 4, 7 \}$.

According to CFL, the corresponding fuzzy sets are the following:
\begin{align*}
C_1 = \begin{array}{|cc|}
x_1;  & w(q_{1, 1}) = 0 \\
x_2;  & w(q_{2, 1}) = 1 \\
x_3;  & w(q_{3, 1}) = 1 \\
x_4;  & w(q_{4, 1}) = 0 \\
x_5;  & w(q_{5, 1}) = 1 \\
x_6;  & w(q_{6, 1}) = 0 \\
x_7;  & w(q_{7, 1}) = 1 \\
\end{array} \hspace{.1in} ; C_2 = \hspace{.1in}  \begin{array}{|cc|}
x_1;  & w(q_{1, 2}) = 1 \\
x_2;  & w(q_{2, 2}) = 1 \\
x_3;  & w(q_{3, 2}) = 0 \\
x_4;  & w(q_{4, 2}) = 1 \\
x_5;  & w(q_{5, 2}) = 0 \\
x_6;  & w(q_{6, 2}) = 0 \\
x_7;  & w(q_{7, 2}) = 1 \\\end{array}\end{align*}

The following operations will be carried out below: a) some operations with classical sets $C_1$ and $C_2$, and b) the same operations with the corresponding sets according to CFL.  In the latter case, this process will be shown only when carrying out these operations with the element $x_7$.\\

\noindent$C_1 \cup C_2 = \langle 2, 3, 5, 7 \rangle \cup \langle 1, 2, 4, 7 \rangle = \langle 1, 2, 3, 4, 5, 7 \rangle$\\ 

The corresponding operation within the CFL framework is as follows:
\begin{align*}
C_1 \cup C_2 = \begin{array}{|cc|}
x_1;  & w(q_{1, 1}) = 0  \\
x_2;  & w(q_{2, 1}) = 1 \\
x_3;  & w(q_{3, 1}) = 1 \\
x_4;  & w(q_{4, 1}) = 0 \\
x_5;  & w(q_{5, 1}) = 1 \\
x_6;  & w(q_{6, 1}) = 0 \\
x_7;  & w(q_{7, 1}) = 1 \\
\end{array} \hspace{.1in} \cup \hspace{.1in}  \begin{array}{|cc|}
x_1;  & w(q_{1, 2}) = 1 \\
x_2;  & w(q_{2, 2}) = 1 \\
x_3;  & w(q_{3, 2}) = 0 \\
x_4;  & w(q_{4, 2}) = 1 \\
x_5;  & w(q_{5, 2}) = 0 \\
x_6;  & w(q_{6, 2}) = 0 \\
x_7;  & w(q_{7, 2}) = 1 \\
\end{array} = \\  \begin{array}{|cc|}
x_1;  & w(q_{1, 1} \lor q_{1, 2}) = 1 \\
x_2;  & w(q_{2, 1} \lor q_{2, 2}) = 1 \\
x_3;  & w(q_{3, 1} \lor q_{3, 2}) = 1 \\
x_4;  & w(q_{4, 1} \lor q_{4, 2}) = 1 \\
x_5;  & w(q_{5, 1} \lor q_{5, 2}) = 1 \\
x_6;  & w(q_{6, 1} \lor q_{6, 2}) = 0 \\
x_7;  & w(q_{7, 1} \lor q_{7, 2}) = 1 \\
\end{array}\end{align*}

\noindent $q_{7, 1}: x_7 \in C_1; w(q_{7,1}) = 1$\\
\noindent $q_{7, 2}: x_7 \in C_2; w(q_{7,2}) = 1$\\

\begin{tabular}{c|c||c|l}
$q_{7, 1}$ & $q_{7, 2}$ & $q_{7, 1} \lor q_{7, 2}$ & $w(q_{7, 1} \lor q_{7, 2}) =\displaystyle\sum_{i=1}^{4}S_i$\\
\midrule
0 & 0 & 0 & $S_1 = 0 $\\
0 & 1 & 1 & $S_2 = w(\overline{q}_{7, 1}) \cdot w(q_{7, 2}) $\\
1 & 0 & 1 & $S_3 = w(q_{7, 1}) \cdot w(\overline{q}_{7, 2}) $\\
1 & 1 & 1 & $S_4 = w(q_{7, 1}) \cdot w(q_{7, 2}) $\\
\end{tabular}\\
\\
\\
\noindent $w(q_{7,1} \lor q_{7,2}) = 0 + (0) \cdot (1) + (1) \cdot (0) + (1) \cdot (1) = 1$\\

\noindent That is, $w(q_{7,1} \lor q_{7,2}) = 1$.\\

\noindent $C_1 \cap C_2 = \langle 2, 3, 5, 7 \rangle \cap \langle 1, 2, 4, 7 \rangle = \langle 2, 7 \rangle$\\

The corresponding operation within the CFL framework is as follows:
\begin{align*}
C_1 \cap C_2 = \begin{array}{|cc|}
x_1;  & w(q_{1, 1}) = 0  \\
x_2;  & w(q_{2, 1}) = 1 \\
x_3;  & w(q_{3, 1}) = 1 \\
x_4;  & w(q_{4, 1}) = 0 \\
x_5;  & w(q_{5, 1}) = 1 \\
x_6;  & w(q_{6, 1}) = 0 \\
x_7;  & w(q_{7, 1}) = 1 \\
\end{array} \hspace{.1in} \cap \hspace{.1in}  \begin{array}{|cc|}
x_1;  & w(q_{1, 2}) = 1 \\
x_2;  & w(q_{2, 2}) = 1 \\
x_3;  & w(q_{3, 2}) = 0 \\
x_4;  & w(q_{4, 2}) = 1 \\
x_5;  & w(q_{5, 2}) = 0 \\
x_6;  & w(q_{6, 2}) = 0 \\
x_7;  & w(q_{7, 2}) = 1 \\
\end{array} = \\  \begin{array}{|cc|}
x_1;  & w(q_{1, 1} \land q_{1, 2}) = 0 \\
x_2;  & w(q_{2, 1} \land q_{2, 2}) = 1 \\
x_3;  & w(q_{3, 1} \land q_{3, 2}) = 0 \\
x_4;  & w(q_{4, 1} \land q_{4, 2}) = 0 \\
x_5;  & w(q_{5, 1} \land q_{5, 2}) = 0 \\
x_6;  & w(q_{6, 1} \land q_{6, 2}) = 0 \\
x_7;  & w(q_{7, 1} \land q_{7, 2}) = 1 \\
\end{array}\end{align*}

\noindent $q_{7, 1}: x_7 \in C_1; w(q_{7,1}) = 1$\\
\noindent $q_{7, 2}: x_7 \in C_2; w(q_{7,2}) = 1$\\

\begin{tabular}{c|c||c|l}
$q_{7, 1}$ & $q_{7, 2}$ & $q_{7, 1} \land q_{7, 2}$ & $w(q_{7, 1} \land q_{7, 2}) =\displaystyle\sum_{i=1}^{4}S_i$\\
\midrule
0 & 0 & 0 & $S_1 = 0 $\\
0 & 1 & 0 & $S_2 = 0$\\
1 & 0 & 0 & $S_3 = 0$\\
1 & 1 & 1 & $S_4 = w(q_{7, 1}) \cdot w(q_{7, 2}) $\\
\end{tabular}\\
\\
\\
\noindent $w(q_{7,1} \land q_{7,2}) = 0 + 0 + 0 + (1) \cdot (1) = 1$\\

\noindent That is, $w(q_{7,1} \land q_{7,2}) = 1$.\\

\noindent $C_1 \dot{\cup} C_2 = \langle 2, 3, 5, 7 \rangle {\dot \cup} \langle 1, 2, 4, 7 \rangle = \langle 1, 3, 4, 5 \rangle$\\ 

The corresponding operation within the CFL framework is as follows:
\begin{align*}
C_1 \dot{\cup} C_2 = \begin{array}{|cc|}
x_1;  & w(q_{1, 1}) = 0  \\
x_2;  & w(q_{2, 1}) = 1 \\
x_3;  & w(q_{3, 1}) = 1 \\
x_4;  & w(q_{4, 1}) = 0 \\
x_5;  & w(q_{5, 1}) = 1 \\
x_6;  & w(q_{6, 1}) = 0 \\
x_7;  & w(q_{7, 1}) = 1 \\
\end{array} \hspace{.1in} \dot{\cup} \hspace{.1in}  \begin{array}{|cc|}
x_1;  & w(q_{1, 2}) = 1 \\
x_2;  & w(q_{2, 2}) = 1 \\
x_3;  & w(q_{3, 2}) = 0 \\
x_4;  & w(q_{4, 2}) = 1 \\
x_5;  & w(q_{5, 2}) = 0 \\
x_6;  & w(q_{6, 2}) = 0 \\
x_7;  & w(q_{7, 2}) = 1 \\
\end{array} = \\  \begin{array}{|cc|}
x_1;  & w(q_{1, 1} \dot{\lor} q_{1, 2}) = 1 \\
x_2;  & w(q_{2, 1} \dot{\lor} q_{2, 2}) = 0 \\
x_3;  & w(q_{3, 1} \dot{\lor} q_{3, 2}) = 1 \\
x_4;  & w(q_{4, 1} \dot{\lor} q_{4, 2}) = 1 \\
x_5;  & w(q_{5, 1} \dot{\lor} q_{5, 2}) = 1 \\
x_6;  & w(q_{6, 1} \dot{\lor} q_{6, 2}) = 0 \\
x_7;  & w(q_{7, 1} \dot{\lor} q_{7, 2}) = 0 \\
\end{array}\end{align*}

\noindent $q_{7, 1}: x_7 \in C_1; w(q_{7,1}) = 1$\\
\noindent $q_{7, 2}: x_7 \in C_2; w(q_{7,2}) = 1$\\

\begin{tabular}{c|c||c|l}
$q_{7, 1}$ & $q_{7, 2}$ & $q_{7, 1} \dot{\lor} q_{7, 2}$ & $w(q_{7, 1} \dot{\lor} q_{7, 2}) =\displaystyle\sum_{i=1}^{4}S_i$\\
\midrule
0 & 0 & 0 & $S_1= 0 $\\
0 & 1 & 1 & $S_2= w(\overline{q}_{7, 1}) \cdot w(q_{7, 2}) $\\
1 & 0 & 1 & $S_3= w(q_{7, 1}) \cdot w(\overline{q}_{7, 2}) $\\
1 & 1 & 0 & $S_4= 0 $\\
\end{tabular}\\
\\

\noindent $w(q_{7,1} \dot{\lor} q_{7,2}) = 0 + (0) \cdot (1) + (1) \cdot (0) + 0 = 0$\\

\noindent That is, $w(q_{7,1} \dot{\lor} q_{7,2}) = 0$.\\
\indent From the semantic viewpoint, there is no difference between any classical set $C$ and the corresponding fuzzy set according to CFL. In effect, suppose first that any element $x$ whatsoever of the universal set $\mathbb U$ considered does belong to $C$. In this case, that element $x$ also belongs to the fuzzy set corresponding to $C$, which will be called $~^fC$, and the proposition associated with $x$ which establishes that membership -- $x \in ~^fC$, abbreviated as $q_{x, ~^fC}$ -- has a weight of truth equal to 1: $w(q_{x, ~^fC}) = 1$. That weight of truth means that $q_{x, ~^fC}$ is ``absolutely'' true.

Second, suppose that the element $x$ does not belong to the classical set $C$. In this case that element $x$ belongs to$~^fC$ such that $w(q_{x, ~^fC}) = 0$. This equality means that the proposition that indicates that $x$ belongs to$~^fC$ (i.e., $q_{x, ~^fC}$) is ``absolutely'' false. That is, here it can also be considered that $x$ does not belong to$~^fC$.

Given that $x$ is any element in the set $\mathbb U$ considered, it can be affirmed that $C = ~^fC$.

The results of operations with classical sets lead to results corresponding to them within the framework of set theory according to CFL. In this sense, it can be considered that the theory of classical sets is a particular case of fuzzy set theory according to CFL.

\section{Discussion and Prospects}

A description was provided of a) how propositional calculus in classical bivalent logic -- BL -- is a particular case of propositional calculus in CFL; and b) how the theory of classical sets is a particular case of set theory in CFL.

Software will be made available to facilitate operations in CFL and in the corresponding fuzzy set theory.

In future articles other operations of fuzzy set theory according to CFL -- different from those already considered -- will be introduced; and a treatment of certain aspects of predicate calculus, according to that variant of fuzzy logic, will be presented. 

This article is the result of an ongoing research program in logic. Its main objective is to establish how BL can be considered, in different ways, a particular case -- or ``limit'' case -- of diverse variants of non-classical logics.

\end{document}